\documentclass[a4paper,twoside,11pt,english,fleqn]{article}

\usepackage[all]{xy}
\usepackage[intlimits]{amsmath}
\usepackage{amssymb,amsthm}
\usepackage{graphicx}
\usepackage[dvips,dvipdf,pdftex]{epsfig}
\usepackage{babel}
\usepackage[T1]{fontenc}
\usepackage{fancyhdr,stmaryrd,color}
\usepackage{times,euler,euscript,eufrak}

  \newtheorem{thm}{Theorem}[section]
  \newtheorem{cor}[thm]{Corollary}
  \newtheorem{lem}[thm]{Lemma}
  \newtheorem{prop}[thm]{Proposition}
  \newtheorem{thm*}{Theorem}[subsection]
  \newtheorem{cor*}[thm*]{Corollary}
  \newtheorem{lem*}[thm*]{Lemma}
  \newtheorem{prop*}[thm*]{Proposition}
\theoremstyle{definition}
  \newtheorem{defn}[thm]{Definition}

  \newtheorem{rem}[thm]{Remark}
  
  \newtheorem{defn*}[thm*]{Definition}
  \newtheorem{ex*}[thm*]{Example}
  \newtheorem{exs*}[thm*]{Examples}
  \newtheorem{rem*}[thm*]{Remark}
  \newtheorem{ntt*}[thm*]{Notation}
\theoremstyle{remark}
  \newtheorem*{dem}{Proof}

\newcommand{\Cr}{\EuScript{C}} 
\newcommand{\Dr}{\EuScript{D}} 
\newcommand{\Or}{\EuScript{O}}

\newcommand{\Fb}{\mathbb{F}} 
\newcommand{\Nb}{\mathbb{N}} 
\newcommand{\Tb}{\mathbb{T}} 
\newcommand{\Vb}{\mathbb{V}} 
\newcommand{\Zb}{\mathbb{Z}}
\newcommand{\zb}{\Zb/2\Zb}
\newcommand{\lz}{\mathrm{L}(\Zb_2)}
\DeclareMathAlphabet{\mathscr}{T1}{pzc}{m}{it}

\newcommand{\red}[1]{\rightarrow\!\!_{{\scriptscriptstyle #1}}}
\newcommand{\mred}[1]{\twoheadrightarrow\!\!_{{\scriptscriptstyle #1}}}
\newcommand{\pmred}[1]{\twoheadrightarrow\!\!_{{\scriptscriptstyle #1}}^+}

\newcommand{\fl}{\rightarrow}
\newcommand{\mfl}{\twoheadrightarrow}

\renewcommand{\phi}{\varphi}
\renewcommand{\epsilon}{\varepsilon}
\newcommand{\tens}{\otimes}
\newcommand{\sq}{\square}
\newcommand{\sur}{\mapsto}
\newcommand{\lsur}{\longmapsto}
\newcommand{\mon}[1]{\langle #1 \rangle}
\newcommand{\ens}[2]{\{#1,\dots,#2\}}
\newcommand{\findem}{\hfill $\diamondsuit$}
\newcommand{\findef}{\hfill $\blacklozenge$}
\newcommand{\ul}[1]{\underline{#1}}
\newdir{ >}{{}*!/-10pt/@{>}}

\newcommand{\emptysectionmark}[1]{\markboth{\textbf{#1}}{\textbf{#1}}}

\fancypagestyle{plain}{
  \fancyhf{}\fancyfoot[LC,RC]{}
  \fancyfoot[LE,RO]{\textbf{\thepage}}
  
  }

\newcommand{\titre}[2]{
\thispagestyle{plain}
\hfill {\large \textbf{#2}} \\
\hrule height 1pt
\medskip
{\huge \textbf{#1}} 
\par\smallskip
\indent {\LARGE \textbf{Yves GUIRAUD}\footnote{Institut de mathématiques de Luminy, Marseille, France - \texttt{guiraud@iml.univ-mrs.fr} }} \\
\hrule height 1pt
}

\begin{document}
\titre{Termination orders for $\mathbf{3}$-dimensional rewriting}{16th November 2004}

\begin{quote}
\textbf{Abstract:} This paper studies $3$-polygraphs as a framework for rewriting on two-dimen-sional words. A translation of term rewriting systems into $3$-polygraphs with explicit resource management is given, and the respective computational properties of each system are studied. Finally, a convergent $3$-polygraph for the (commutative) theory of $\zb$-vector spaces is given. In order to prove these results, it is explained how to craft a class of termination orders for $3$-polygraphs.
\end{quote}

\section*{Outline}
\emptysectionmark{Outline}

This paper starts with the introductory section~\ref{equational} on equational theories and term rewriting systems. It gives notations and graphical representations that are used in the sequel. Then, it focuses on one major restriction of term rewriting, namely the fact that it cannot provide convergent presentations for \emph{commutative} equational theories: equational theories that contain a commutative binary operator.

Section~\ref{resource} studies the resource management operations of permutation, erasure and duplication: they are implicit and global in term rewriting and it is sketched there how to make them explicit. However, the framework for rewriting in algebraic structures needs to be extended to include this change; section~\ref{operads} proposes $3$-polygraphs to fulfill this role. Here, these objects, introduced in [Burroni 1993], are used as equational presentations of a special case of $2$-categories: MacLane's product categories, called PROs, for short, in [MacLane 1965].

These first three sections do not introduce new material, but focus on the notations, representations, terminology and philosophy of this paper. Then section~\ref{explicit} gives some relations between term rewriting systems and $3$-polygraphs: a translation from the former to the latter is built and some properties are given. The main result of the section is the proof of a conjecture from [Lafont 2003]: any left-linear convergent term rewriting system can be translated into a convergent $3$-polygraph.

To prove some of these results, one needs new tools, in adequation with the more complicated structure of polygraphs. In particular, section~\ref{technique} introduces a recipe to build termination orders for them. Section~\ref{app1} consists in the application of this technique to prove some termination results of section~\ref{explicit}. Finally, section~\ref{app2} applies the same technique to prove the termination of the $3$-polygraph $\lz$ which was introduced in [Lafont 2003] and, since then, was already known to be a confluent presentation of the equational theory of $\zb$-vector spaces. It is therefore the first known convergent presentation of a commutative equational theory.

\section{Equational theories and term rewriting systems}\label{equational}

Universal algebra provides different types of objects in order to modelize algebraic structures. Among them are \emph{equational theories}: these are presentations by generators (or \emph{operators}) and relations (or \emph{equations}, \emph{equalities}). As an example, the equational theory of \emph{monoids} is a pair~$(\Sigma,E_0)$ consisting of the \emph{signature} $\Sigma$ (a set of operators) and the family $E_0$ of equations given by:
\begin{align*}
&\Sigma=\{\:\mu:2\fl 1,\:\eta:0\fl 1\:\}\:, \\
&E_0=\big(\:\mu(\mu(x,y),z)=\mu(x,\mu(y,z)), \:\mu(\eta,x)=x, \: \mu(x,\eta)=x\:\big).
\end{align*}

\noindent Each operator has a finite number of inputs and of outputs. When each one has exactly one output, which is the case here, the signature is said to be \emph{algebraic}. The given equational theory~$(\Sigma,E_0)$ is said to be \emph{the theory of monoids} since monoids are exactly sets endowed with a binary operation and a constant, such that the operation is associative and admits the constant as a left and right unit. \\

\noindent The formal operations one can form on any set with a binary operation and a constant are called the \emph{terms} built from the signature $\Sigma$. There exist numerous ways to build the set $T\Sigma$ of such terms, and each one gives a different representation for them. Two are used here, a \emph{syntactic} one and a \emph{diagrammatic} one. For each one, a fixed countable set $V$ is needed; its elements are called \emph{variables}.

The classical representation of terms define them inductively with the following construction rules: the first one states that each variable is a term; furthermore, the constant $\eta$ is a term; then, for any two terms $u$ and $v$, the formal expression $\mu(u,v)$ is a term.

The diagrammatic representation starts with the assignment, for each operator with $n$ inputs, of an arbitrarily chosen tree of height one with $n$ leaves. For example, one can fix the following trees:
\begin{center}
\begin{picture}(0,0)%
\includegraphics{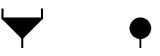}%
\end{picture}%
\setlength{\unitlength}{4144sp}%
\begingroup\makeatletter\ifx\SetFigFont\undefined%
\gdef\SetFigFont#1#2#3#4#5{%
  \reset@font\fontsize{#1}{#2pt}%
  \fontfamily{#3}\fontseries{#4}\fontshape{#5}%
  \selectfont}%
\fi\endgroup%
\begin{picture}(695,421)(169,250)
\put(271,299){\makebox(0,0)[b]{\smash{{\SetFigFont{12}{14.4}{\familydefault}{\mddefault}{\updefault}{\color[rgb]{0,0,0}$\mu$}%
}}}}
\put(811,299){\makebox(0,0)[b]{\smash{{\SetFigFont{12}{14.4}{\familydefault}{\mddefault}{\updefault}{\color[rgb]{0,0,0}$\eta$}%
}}}}
\end{picture}%
\end{center}

\noindent Then, the terms are all the trees one can build from these two generating trees and which leaves are labelled with variables. As an example, the following figure pictures terms built from the signature $\Sigma$, with the two representations for each one:
\begin{center}
\begin{picture}(0,0)%
\includegraphics{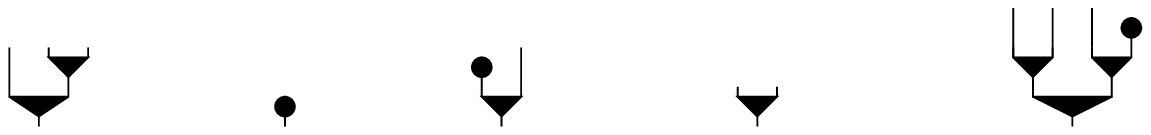}%
\end{picture}%
\setlength{\unitlength}{4144sp}%
\begingroup\makeatletter\ifx\SetFigFont\undefined%
\gdef\SetFigFont#1#2#3#4#5{%
  \reset@font\fontsize{#1}{#2pt}%
  \fontfamily{#3}\fontseries{#4}\fontshape{#5}%
  \selectfont}%
\fi\endgroup%
\begin{picture}(5322,1069)(312,-209)
\put(5311,-151){\makebox(0,0)[b]{\smash{{\SetFigFont{12}{14.4}{\familydefault}{\mddefault}{\updefault}{\color[rgb]{0,0,0}$\mu(\mu(x,y),\mu(x,\eta))$}%
}}}}
\put(451,524){\makebox(0,0)[b]{\smash{{\SetFigFont{12}{14.4}{\familydefault}{\mddefault}{\updefault}{\color[rgb]{0,0,0}$x$}%
}}}}
\put(631,524){\makebox(0,0)[b]{\smash{{\SetFigFont{12}{14.4}{\familydefault}{\mddefault}{\updefault}{\color[rgb]{0,0,0}$y$}%
}}}}
\put(811,524){\makebox(0,0)[b]{\smash{{\SetFigFont{12}{14.4}{\familydefault}{\mddefault}{\updefault}{\color[rgb]{0,0,0}$z$}%
}}}}
\put(2791,524){\makebox(0,0)[b]{\smash{{\SetFigFont{12}{14.4}{\familydefault}{\mddefault}{\updefault}{\color[rgb]{0,0,0}$x$}%
}}}}
\put(3781,344){\makebox(0,0)[b]{\smash{{\SetFigFont{12}{14.4}{\familydefault}{\mddefault}{\updefault}{\color[rgb]{0,0,0}$x$}%
}}}}
\put(3961,344){\makebox(0,0)[b]{\smash{{\SetFigFont{12}{14.4}{\familydefault}{\mddefault}{\updefault}{\color[rgb]{0,0,0}$x$}%
}}}}
\put(5041,704){\makebox(0,0)[b]{\smash{{\SetFigFont{12}{14.4}{\familydefault}{\mddefault}{\updefault}{\color[rgb]{0,0,0}$x$}%
}}}}
\put(5401,704){\makebox(0,0)[b]{\smash{{\SetFigFont{12}{14.4}{\familydefault}{\mddefault}{\updefault}{\color[rgb]{0,0,0}$x$}%
}}}}
\put(5221,704){\makebox(0,0)[b]{\smash{{\SetFigFont{12}{14.4}{\familydefault}{\mddefault}{\updefault}{\color[rgb]{0,0,0}$y$}%
}}}}
\put(586,-151){\makebox(0,0)[b]{\smash{{\SetFigFont{12}{14.4}{\familydefault}{\mddefault}{\updefault}{\color[rgb]{0,0,0}$\mu(x,\mu(y,z))$}%
}}}}
\put(1711,-151){\makebox(0,0)[b]{\smash{{\SetFigFont{12}{14.4}{\familydefault}{\mddefault}{\updefault}{\color[rgb]{0,0,0}$\eta$}%
}}}}
\put(2701,-151){\makebox(0,0)[b]{\smash{{\SetFigFont{12}{14.4}{\familydefault}{\mddefault}{\updefault}{\color[rgb]{0,0,0}$\mu(\eta,x)$}%
}}}}
\put(3871,-151){\makebox(0,0)[b]{\smash{{\SetFigFont{12}{14.4}{\familydefault}{\mddefault}{\updefault}{\color[rgb]{0,0,0}$\mu(x,x)$}%
}}}}
\end{picture}%
\end{center}

\noindent The equations from the theory of monoids generate equalities between terms that represent the same operation, through a \emph{rewriting} process. Let us sketch how this works. For example, the following term contains the tree-part of the associativity rule left-member, which has been greyed out:
\begin{center}
\begin{picture}(0,0)%
\includegraphics{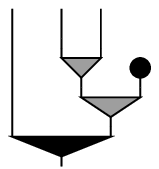}%
\end{picture}%
\setlength{\unitlength}{4144sp}%
\begingroup\makeatletter\ifx\SetFigFont\undefined%
\gdef\SetFigFont#1#2#3#4#5{%
  \reset@font\fontsize{#1}{#2pt}%
  \fontfamily{#3}\fontseries{#4}\fontshape{#5}%
  \selectfont}%
\fi\endgroup%
\begin{picture}(1000,933)(-319,17)
\put(-89,794){\makebox(0,0)[b]{\smash{{\SetFigFont{12}{14.4}{\familydefault}{\mddefault}{\updefault}{\color[rgb]{0,0,0}$x_1$}%
}}}}
\put(181,794){\makebox(0,0)[b]{\smash{{\SetFigFont{12}{14.4}{\familydefault}{\mddefault}{\updefault}{\color[rgb]{0,0,0}$x_2$}%
}}}}
\put(451,794){\makebox(0,0)[b]{\smash{{\SetFigFont{12}{14.4}{\familydefault}{\mddefault}{\updefault}{\color[rgb]{0,0,0}$x_3$}%
}}}}
\end{picture}%
\end{center}

\noindent Hence, the associativity equation generates an equality between the chosen term and another. To determine which one, let us follow the following method, which consists of three steps: at first, the remaining (black) part of the term is copied; then, in the space left empty, the other member of the rule is placed; finally, the two parts obtained are joined (by dotted lines), according to the respective position of the variables in each member of the equation. Concerning our example, this process is pictured as follows:
\begin{center}
\begin{picture}(0,0)%
\includegraphics{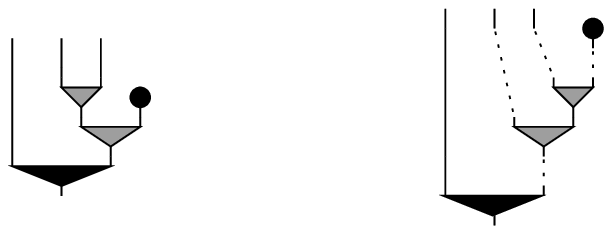}%
\end{picture}%
\setlength{\unitlength}{4144sp}%
\begingroup\makeatletter\ifx\SetFigFont\undefined%
\gdef\SetFigFont#1#2#3#4#5{%
  \reset@font\fontsize{#1}{#2pt}%
  \fontfamily{#3}\fontseries{#4}\fontshape{#5}%
  \selectfont}%
\fi\endgroup%
\begin{picture}(2983,1203)(-319,-118)
\put(1891,929){\makebox(0,0)[b]{\smash{{\SetFigFont{12}{14.4}{\familydefault}{\mddefault}{\updefault}{\color[rgb]{0,0,0}$x_1$}%
}}}}
\put(181,794){\makebox(0,0)[b]{\smash{{\SetFigFont{12}{14.4}{\familydefault}{\mddefault}{\updefault}{\color[rgb]{0,0,0}$x_2$}%
}}}}
\put(451,794){\makebox(0,0)[b]{\smash{{\SetFigFont{12}{14.4}{\familydefault}{\mddefault}{\updefault}{\color[rgb]{0,0,0}$x_3$}%
}}}}
\put(-89,794){\makebox(0,0)[b]{\smash{{\SetFigFont{12}{14.4}{\familydefault}{\mddefault}{\updefault}{\color[rgb]{0,0,0}$x_1$}%
}}}}
\put(2161,929){\makebox(0,0)[b]{\smash{{\SetFigFont{12}{14.4}{\familydefault}{\mddefault}{\updefault}{\color[rgb]{0,0,0}$x_2$}%
}}}}
\put(2431,929){\makebox(0,0)[b]{\smash{{\SetFigFont{12}{14.4}{\familydefault}{\mddefault}{\updefault}{\color[rgb]{0,0,0}$x_3$}%
}}}}
\end{picture}%
\end{center}

\noindent Note that each variable appears once and in the same position in each member of the associativity rule, so that the links are direct. When the second term is compacted, the following equality holds and is said to be generated by the associativity equation:
\begin{center}
\begin{picture}(0,0)%
\includegraphics{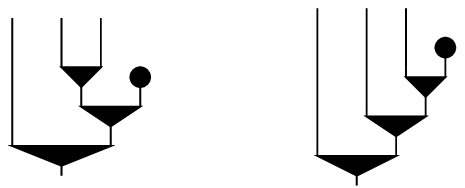}%
\end{picture}%
\setlength{\unitlength}{4144sp}%
\begingroup\makeatletter\ifx\SetFigFont\undefined%
\gdef\SetFigFont#1#2#3#4#5{%
  \reset@font\fontsize{#1}{#2pt}%
  \fontfamily{#3}\fontseries{#4}\fontshape{#5}%
  \selectfont}%
\fi\endgroup%
\begin{picture}(2395,1068)(-319,-28)
\put(1261,884){\makebox(0,0)[b]{\smash{{\SetFigFont{12}{14.4}{\familydefault}{\mddefault}{\updefault}{\color[rgb]{0,0,0}$x_1$}%
}}}}
\put(181,794){\makebox(0,0)[b]{\smash{{\SetFigFont{12}{14.4}{\rmdefault}{\mddefault}{\updefault}{\color[rgb]{0,0,0}$x_2$}%
}}}}
\put(451,794){\makebox(0,0)[b]{\smash{{\SetFigFont{12}{14.4}{\rmdefault}{\mddefault}{\updefault}{\color[rgb]{0,0,0}$x_3$}%
}}}}
\put(-89,794){\makebox(0,0)[b]{\smash{{\SetFigFont{12}{14.4}{\rmdefault}{\mddefault}{\updefault}{\color[rgb]{0,0,0}$x_1$}%
}}}}
\put(946,344){\makebox(0,0)[b]{\smash{{\SetFigFont{12}{14.4}{\rmdefault}{\mddefault}{\updefault}{\color[rgb]{0,0,0}$=$}%
}}}}
\put(1576,884){\makebox(0,0)[b]{\smash{{\SetFigFont{12}{14.4}{\familydefault}{\mddefault}{\updefault}{\color[rgb]{0,0,0}$x_2$}%
}}}}
\put(1846,884){\makebox(0,0)[b]{\smash{{\SetFigFont{12}{14.4}{\familydefault}{\mddefault}{\updefault}{\color[rgb]{0,0,0}$x_3$}%
}}}}
\end{picture}%
\end{center}

\noindent In order to study the computational properties of these rewriting processes, \emph{term rewriting systems} are useful; they can be defined as oriented equational theories. Indeed, such a rewriting system is defined from an equational theory by keeping the same operators and replacing each equation by a \emph{rewrite rule}: it is an oriented version of the equation, which can only be used in one way. As an example, starting from the equational theory of monoids, one can form the term rewriting system $(\Sigma,R_0)$, where $\Sigma$ is still the same algebraic signature made of a product $\mu$ and a unit $\eta$ and $R_0$ is the following set of three rules:
$$
\mu(\mu(x,y),z)\fl\mu(x,\mu(y,z)), \quad \mu(\eta,x)\fl x, \quad \mu(x,\eta)\fl x.
$$

\noindent Rewrite rules generate \emph{reductions} instead of equalities, and a graph containing terms as vertices and reductions as edges is called a \emph{reduction graph}. Some geometrical properties of reduction graphs are of particular interest since they have consequences on computational properties of the rewriting process. Among these geometrical properties, three are particularly studied: \emph{termination}, \emph{confluence} and \emph{convergence}. \\

\noindent A rewriting system \emph{terminates} if it contains no infinite length reduction paths such as:
$$
u_0\fl u_1\fl u_2\fl\dots\fl u_n\fl u_{n+1}\fl \dots
$$

\noindent Intuitively, this means that the rewriting calculus must end after a finite time, whatever the input is. This is formalized by the following consequence of termination: every term $u$ has at least one \emph{normal form} $\hat{u}$; this means that $\hat{u}$ is a term such that there exists a finite reduction path from $u$ to $\hat{u}$ (denoted by $u\mfl\hat{u}$) and $\hat{u}$ is irreducible (no rule can apply on it). \\

\noindent A rewriting system is \emph{confluent} if, whenever there exist three terms $u$, $v$ and $w$ such that~$u\mfl v$ and~$u\mfl w$, then there exists a fourth term $t$ such that $v\mfl t$ and $w\mfl t$. Intuitively, this means that choices made between two rules that can transform the same term do not have any consequence on a potential final result; equivalently, this means that any term has at most one normal form. \\

\noindent Thus, one defines the last property: a rewriting system is \emph{convergent} when it is both terminating and confluent. One immediate consequence is that any term has exactly one normal form. This property is very useful for several purposes.

One of the most known is the following usage: let us assume that $(\Sigma,E)$ is an equational theory and that $(\Sigma,R)$ is a rewriting system that is a \emph{finite convergent presentation} of $(\Sigma,R)$, which means that it is a convergent rewriting system with a finite number of rules and such that two terms are equal in the equational theory if and only if there exists a non oriented reduction path between these two terms in the rewriting system. Then there exists a decision procedure to check if two terms $u$ and $v$ are equal or not.

Indeed, one computes their unique normal forms $\hat{u}$ and $\hat{v}$. Note that this is where the finiteness condition is useful: it allows one to check if a term is a normal form. Then the two normal forms $\hat{u}$ and $\hat{v}$ are compared: $u$ and $v$ are equal in the equational theory if and only if $\hat{u}$ and $\hat{v}$ are (synctactically) equal. \\

\noindent However, term rewriting systems have a major restriction in this field: there is a large class of equational theories for which they cannot provide a convergent presentation. These are the commutative theories, fairly frequent in algebra, which are equational theories with a commutative binary operator. As an example, let us take a look at one of the simplest, namely the equational theory of commutative monoids. Its signature is still~$\Sigma$; its set $E_1$ of equations is made of the same three as the ones for monoids (associativity and left and right units) plus the following one expressing the commutativity of the product:
$$
\mu(x,y)=\mu(y,x).
$$

\noindent From this theory, one can form a number of term rewriting systems, such as the one with $\Sigma$ as signature and with the following choice $R_1$ of orientations for equations:
$$
\mu(\mu(x,y),z)\fl\mu(x,\mu(y,z)), \quad \mu(\eta,x)\fl x, \quad \mu(x,\eta)\fl x, \quad \mu(x,y)\fl\mu(y,x).
$$

\noindent Note that the last rule could have been chosen in the reverse direction, but it would not change the following fact: this rule generates infinite reduction paths. Indeed, for any two terms $u$ and~$v$, the commutativity rules generates:
$$
\mu(u,v)\fl\mu(v,u)\fl\mu(u,v)\fl\mu(v,u)\fl\dots
$$

\noindent The purpose of this paper is to provide a framework where some commutative equational theories admit convergent presentations: $3$-polygraphs. Links between term rewriting systems and $3$-polygraphs are studied and a new tool to prove termination is given and applied on some examples. \\

\noindent The equational theory that provides the main example here is the one of $\zb$-vector spaces: it has the same operators as the previous ones (the binary product embodies the sum and the unit is the zero) and a set $E_2$ of five equations made of the four from $E_1$ (associativity, left and right units and commutativity) plus the following fifth equation:
$$
\mu(x,x)=\eta.
$$

\noindent It expresses the fact that, in a $\zb$-vector space, any element is its own opposite. This theory is prefered to the theory of commutative monoids for two reasons. The first one is theoretical: any boolean algebra has an underlying $\zb$-vector space, so that any convergent presentation for $\zb$-vector spaces is a first step towards one for boolean circuits. The second one concerns the application range of the tools developped here: this fifth equation has some nasty computational effects and is thus important to encompass in the new framework, so that it can be used for other applications. \\

\noindent From the theory of $\zb$-vector spaces, the term rewriting system $(\Sigma,R_2)$ is built, where $R_2$ is the following choice of orientations:
$$
\mu(\mu(x,y),z)\fl\mu(x,\mu(y,z)), \quad \mu(\eta,x)\fl x, \quad \mu(x,\eta)\fl x, \quad \mu(x,y)\fl\mu(y,x), \quad \mu(x,x)\fl\eta.
$$

\noindent Note that this rewriting system is neither terminating nor confluent but will serve as a starting point to build a convergent presentation. This transformation will start with the study of the so-called \emph{resource management operations}. For further information on (term) rewriting systems, one can refer to [Baader Nipkow 1998].

\section{Resource management operations}\label{resource}

\noindent Let us recall the last step of the term rewriting process: one has to draw links between two parts of a term, according to the variables occuring in the corresponding rule. As mentionned earlier, the rewriting example in section \ref{equational} is the simpliest case: indeed, the variables occur once each and in the same order in each member of the associativity rule. However, if this is not the case, one has two use additional operations before links are drawn: these operations are called the \emph{resource management operations} and there are three of the kind, \emph{permutation}, \emph{erasure} and \emph{duplication}. \\

\noindent Permutation is used, for example, when the commutativity rule is applied. Indeed, when in this case, one has to use a permutation operation that will exchange the two grey subterms in any term such as the following generic one:
\begin{center}
\begin{picture}(0,0)%
\includegraphics{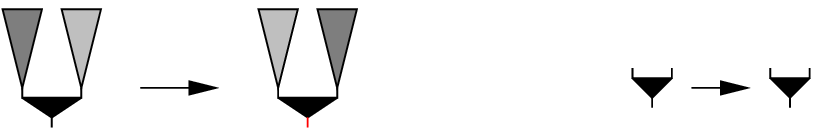}%
\end{picture}%
\setlength{\unitlength}{4144sp}%
\begingroup\makeatletter\ifx\SetFigFont\undefined%
\gdef\SetFigFont#1#2#3#4#5{%
  \reset@font\fontsize{#1}{#2pt}%
  \fontfamily{#3}\fontseries{#4}\fontshape{#5}%
  \selectfont}%
\fi\endgroup%
\begin{picture}(3840,564)(169,17)
\put(3061,344){\makebox(0,0)[b]{\smash{{\SetFigFont{12}{14.4}{\familydefault}{\mddefault}{\updefault}{\color[rgb]{0,0,0}$x$}%
}}}}
\put(3871,344){\makebox(0,0)[b]{\smash{{\SetFigFont{12}{14.4}{\familydefault}{\mddefault}{\updefault}{\color[rgb]{0,0,0}$x$}%
}}}}
\put(3241,344){\makebox(0,0)[b]{\smash{{\SetFigFont{12}{14.4}{\familydefault}{\mddefault}{\updefault}{\color[rgb]{0,0,0}$y$}%
}}}}
\put(3691,344){\makebox(0,0)[b]{\smash{{\SetFigFont{12}{14.4}{\familydefault}{\mddefault}{\updefault}{\color[rgb]{0,0,0}$y$}%
}}}}
\end{picture}%
\end{center}

\noindent The second operation, erasure, is used in the following case, for example: let us consider a theory containing a binary operator and a constant which is a right absorbing element. The following figure displays a rule which expresses this property (on the right) together with a generic application of this rule (on the left); this requires an intermediate operation that erases the grey subterm:
\begin{center}
\begin{picture}(0,0)%
\includegraphics{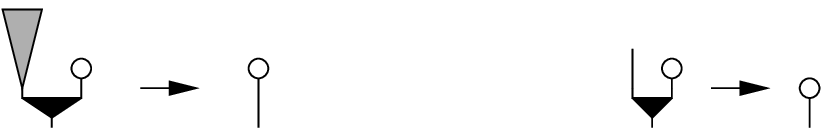}%
\end{picture}%
\setlength{\unitlength}{4144sp}%
\begingroup\makeatletter\ifx\SetFigFont\undefined%
\gdef\SetFigFont#1#2#3#4#5{%
  \reset@font\fontsize{#1}{#2pt}%
  \fontfamily{#3}\fontseries{#4}\fontshape{#5}%
  \selectfont}%
\fi\endgroup%
\begin{picture}(3755,564)(169,-73)
\put(3061,344){\makebox(0,0)[b]{\smash{{\SetFigFont{12}{14.4}{\familydefault}{\mddefault}{\updefault}{\color[rgb]{0,0,0}$x$}%
}}}}
\end{picture}%
\end{center}

\noindent Finally, the last operation, called duplication, can occur in the following case: let us consider a theory containing two binary operators, one of which is left-distributive with respect to the other. Then, when applied, a rule that expresses this property (such as the one pictured on the right) requires the use of an operation that can duplicate the greymost subterm (and exchange one of its copies with another subterm, but this is the already-encountered permutation):
\begin{center}
\begin{picture}(0,0)%
\includegraphics{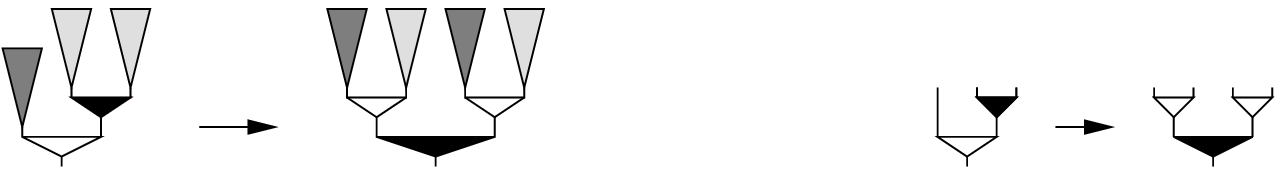}%
\end{picture}%
\setlength{\unitlength}{4144sp}%
\begingroup\makeatletter\ifx\SetFigFont\undefined%
\gdef\SetFigFont#1#2#3#4#5{%
  \reset@font\fontsize{#1}{#2pt}%
  \fontfamily{#3}\fontseries{#4}\fontshape{#5}%
  \selectfont}%
\fi\endgroup%
\begin{picture}(5955,744)(304,-163)
\put(2701,434){\makebox(0,0)[b]{\smash{{\SetFigFont{10}{12.0}{\familydefault}{\mddefault}{\updefault}{\color[rgb]{0,0,0}$2$}%
}}}}
\put(4591,254){\makebox(0,0)[b]{\smash{{\SetFigFont{12}{14.4}{\familydefault}{\mddefault}{\updefault}{\color[rgb]{0,0,0}$x$}%
}}}}
\put(5581,254){\makebox(0,0)[b]{\smash{{\SetFigFont{12}{14.4}{\familydefault}{\mddefault}{\updefault}{\color[rgb]{0,0,0}$x$}%
}}}}
\put(5941,254){\makebox(0,0)[b]{\smash{{\SetFigFont{12}{14.4}{\familydefault}{\mddefault}{\updefault}{\color[rgb]{0,0,0}$x$}%
}}}}
\put(4771,254){\makebox(0,0)[b]{\smash{{\SetFigFont{12}{14.4}{\familydefault}{\mddefault}{\updefault}{\color[rgb]{0,0,0}$y$}%
}}}}
\put(5761,254){\makebox(0,0)[b]{\smash{{\SetFigFont{12}{14.4}{\familydefault}{\mddefault}{\updefault}{\color[rgb]{0,0,0}$y$}%
}}}}
\put(4951,254){\makebox(0,0)[b]{\smash{{\SetFigFont{12}{14.4}{\familydefault}{\mddefault}{\updefault}{\color[rgb]{0,0,0}$z$}%
}}}}
\put(6121,254){\makebox(0,0)[b]{\smash{{\SetFigFont{12}{14.4}{\familydefault}{\mddefault}{\updefault}{\color[rgb]{0,0,0}$z$}%
}}}}
\put(631,434){\makebox(0,0)[b]{\smash{{\SetFigFont{10}{12.0}{\familydefault}{\mddefault}{\updefault}{\color[rgb]{0,0,0}$1$}%
}}}}
\put(901,434){\makebox(0,0)[b]{\smash{{\SetFigFont{10}{12.0}{\familydefault}{\mddefault}{\updefault}{\color[rgb]{0,0,0}$2$}%
}}}}
\put(2161,434){\makebox(0,0)[b]{\smash{{\SetFigFont{10}{12.0}{\familydefault}{\mddefault}{\updefault}{\color[rgb]{0,0,0}$1$}%
}}}}
\end{picture}%
\end{center}

\noindent Thus, in term rewriting, these three operations are both \emph{implicit} (they are not specified by rules) and \emph{global} (they act immediately on subterms of any size). We are now going to sketch how one can make them explicit and local: only the idea is given here, the full translation is postponed to section \ref{explicit}. \\

\noindent Let us start with the following observation: the use of the three resource management operations is specified both by the number of occurences and the order of appearance of each variable in each member of a rewrite rule. Thus, in order to make these operations explicit, variables will be replaced by some additional operators that will represent local permutations, erasers and duplicators; furthermore, rules will guarantee the global behaviour of these local operators. \\

\noindent In order to give an idea of how the translation works, let us start with the study of this term, which represents the operation $(x,y,z)\sur\mu(\mu(x,z),x)$:
\begin{center}
\begin{picture}(0,0)%
\includegraphics{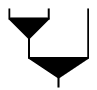}%
\end{picture}%
\setlength{\unitlength}{4144sp}%
\begingroup\makeatletter\ifx\SetFigFont\undefined%
\gdef\SetFigFont#1#2#3#4#5{%
  \reset@font\fontsize{#1}{#2pt}%
  \fontfamily{#3}\fontseries{#4}\fontshape{#5}%
  \selectfont}%
\fi\endgroup%
\begin{picture}(637,549)(42,-523)
\put(541,-106){\makebox(0,0)[b]{\smash{{\SetFigFont{12}{14.4}{\familydefault}{\mddefault}{\updefault}{\color[rgb]{0,0,0}$x$}%
}}}}
\put(181,-106){\makebox(0,0)[b]{\smash{{\SetFigFont{12}{14.4}{\familydefault}{\mddefault}{\updefault}{\color[rgb]{0,0,0}$x$}%
}}}}
\put(361,-106){\makebox(0,0)[b]{\smash{{\SetFigFont{12}{14.4}{\familydefault}{\mddefault}{\updefault}{\color[rgb]{0,0,0}$z$}%
}}}}
\end{picture}%
\end{center}

\noindent Seen as an operation, it is the composite of $(x,y,z)\sur(x,z,x)$ followed by $(x,y,z)\sur\mu(\mu(x,y),z)$. The first operation can be pictured as the following diagram (a shunter), since its action is to tell where each of the three arguments goes in the term:
\begin{center}
\begin{picture}(0,0)%
\includegraphics{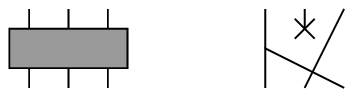}%
\end{picture}%
\setlength{\unitlength}{4144sp}%
\begingroup\makeatletter\ifx\SetFigFont\undefined%
\gdef\SetFigFont#1#2#3#4#5{%
  \reset@font\fontsize{#1}{#2pt}%
  \fontfamily{#3}\fontseries{#4}\fontshape{#5}%
  \selectfont}%
\fi\endgroup%
\begin{picture}(1716,745)(-47,-155)
\put(1351,-106){\makebox(0,0)[b]{\smash{{\SetFigFont{12}{14.4}{\familydefault}{\mddefault}{\updefault}{\color[rgb]{0,0,0}$z$}%
}}}}
\put( 91,434){\makebox(0,0)[b]{\smash{{\SetFigFont{12}{14.4}{\familydefault}{\mddefault}{\updefault}{\color[rgb]{0,0,0}$x$}%
}}}}
\put(1171,434){\makebox(0,0)[b]{\smash{{\SetFigFont{12}{14.4}{\familydefault}{\mddefault}{\updefault}{\color[rgb]{0,0,0}$x$}%
}}}}
\put( 91,-106){\makebox(0,0)[b]{\smash{{\SetFigFont{12}{14.4}{\familydefault}{\mddefault}{\updefault}{\color[rgb]{0,0,0}$x$}%
}}}}
\put(451,-106){\makebox(0,0)[b]{\smash{{\SetFigFont{12}{14.4}{\familydefault}{\mddefault}{\updefault}{\color[rgb]{0,0,0}$x$}%
}}}}
\put(1171,-106){\makebox(0,0)[b]{\smash{{\SetFigFont{12}{14.4}{\familydefault}{\mddefault}{\updefault}{\color[rgb]{0,0,0}$x$}%
}}}}
\put(1531,-106){\makebox(0,0)[b]{\smash{{\SetFigFont{12}{14.4}{\familydefault}{\mddefault}{\updefault}{\color[rgb]{0,0,0}$x$}%
}}}}
\put(271,434){\makebox(0,0)[b]{\smash{{\SetFigFont{12}{14.4}{\familydefault}{\mddefault}{\updefault}{\color[rgb]{0,0,0}$y$}%
}}}}
\put(1351,434){\makebox(0,0)[b]{\smash{{\SetFigFont{12}{14.4}{\familydefault}{\mddefault}{\updefault}{\color[rgb]{0,0,0}$y$}%
}}}}
\put(451,434){\makebox(0,0)[b]{\smash{{\SetFigFont{12}{14.4}{\familydefault}{\mddefault}{\updefault}{\color[rgb]{0,0,0}$z$}%
}}}}
\put(271,-106){\makebox(0,0)[b]{\smash{{\SetFigFont{12}{14.4}{\familydefault}{\mddefault}{\updefault}{\color[rgb]{0,0,0}$z$}%
}}}}
\put(1531,434){\makebox(0,0)[b]{\smash{{\SetFigFont{12}{14.4}{\familydefault}{\mddefault}{\updefault}{\color[rgb]{0,0,0}$z$}%
}}}}
\put(856,164){\makebox(0,0)[b]{\smash{{\SetFigFont{12}{14.4}{\familydefault}{\mddefault}{\updefault}{\color[rgb]{0,0,0}$=$}%
}}}}
\end{picture}%
\end{center}

\noindent This diagram will be formalized as a composite of new operators and the term will be translated this way (with some explanations below):
\begin{center}
\begin{picture}(0,0)%
\includegraphics{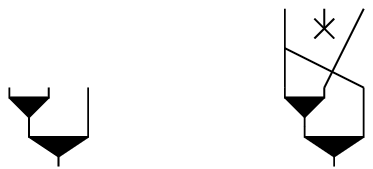}%
\end{picture}%
\setlength{\unitlength}{4144sp}%
\begingroup\makeatletter\ifx\SetFigFont\undefined%
\gdef\SetFigFont#1#2#3#4#5{%
  \reset@font\fontsize{#1}{#2pt}%
  \fontfamily{#3}\fontseries{#4}\fontshape{#5}%
  \selectfont}%
\fi\endgroup%
\begin{picture}(1771,744)(42,-523)
\put(361,-106){\makebox(0,0)[b]{\smash{{\SetFigFont{12}{14.4}{\familydefault}{\mddefault}{\updefault}{\color[rgb]{0,0,0}$3$}%
}}}}
\put(991,-331){\makebox(0,0)[b]{\smash{{\SetFigFont{12}{14.4}{\familydefault}{\mddefault}{\updefault}{\color[rgb]{0,0,0}$\lsur$}%
}}}}
\put(181,-106){\makebox(0,0)[b]{\smash{{\SetFigFont{12}{14.4}{\familydefault}{\mddefault}{\updefault}{\color[rgb]{0,0,0}$1$}%
}}}}
\put(541,-106){\makebox(0,0)[b]{\smash{{\SetFigFont{12}{14.4}{\familydefault}{\mddefault}{\updefault}{\color[rgb]{0,0,0}$1$}%
}}}}
\end{picture}%
\end{center}

\noindent Variables in the term have been replaced by ordinals; indeed, we have seen that variables are just labels corresponding to the first, second, third, etc. arguments taken by the corresponding operation. Hence, they will be replaced by ordinals whenever it makes the translation clearer. The second remark is also about variables, but in the translated diagram: they will always appear, after translation, in order: 1, 2, 3, etc. Thus, they have no purpose anymore; they will therefore vanish, as in the diagram. \\

\noindent Finally, let us see what operators will be added to the signature and sketch how to translate terms and rules. One operator is added for each resource management operation: indeed, in order to formalize our previous diagram, one must be able to exchange two arguments, erase one or duplicate another one. Thus, we fix a (non-algebraic) signature $\Delta$ made of the following three \emph{resource management} operators:
\begin{center}
\begin{picture}(0,0)%
\includegraphics{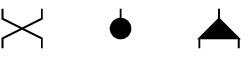}%
\end{picture}%
\setlength{\unitlength}{4144sp}%
\begingroup\makeatletter\ifx\SetFigFont\undefined%
\gdef\SetFigFont#1#2#3#4#5{%
  \reset@font\fontsize{#1}{#2pt}%
  \fontfamily{#3}\fontseries{#4}\fontshape{#5}%
  \selectfont}%
\fi\endgroup%
\begin{picture}(1104,475)(529,16)
\put(1531, 74){\makebox(0,0)[b]{\smash{{\SetFigFont{12}{14.4}{\familydefault}{\mddefault}{\updefault}$\delta$}}}}
\put(631, 74){\makebox(0,0)[b]{\smash{{\SetFigFont{12}{14.4}{\familydefault}{\mddefault}{\updefault}$\tau$}}}}
\put(1081, 74){\makebox(0,0)[b]{\smash{{\SetFigFont{12}{14.4}{\familydefault}{\mddefault}{\updefault}$\epsilon$}}}}
\end{picture}%
\end{center}

\noindent Each one has a representation that makes explicit the operation one wishes it to embody. Some rules will be added to ensure their global behaviour, but they will be given in section \ref{explicit}. For the moment, the only thing we need to know is that these rules give the following interpretations to these three operators:
$$
\tau(x,y)=(y,x), \quad \epsilon(x)=(\text{\texttt{nothing}}), \quad \delta(x)=(x,x).
$$

\noindent Now, let us sketch how terms are translated: first, the tree-part is copied; then and progressively, resource management operators are added on the top of the copy, according to the variables that appear in the term. The following figure gives four sample translations (the translating map is denoted by~$\Phi$ thereafter):
\begin{center}
\begin{picture}(0,0)%
\includegraphics{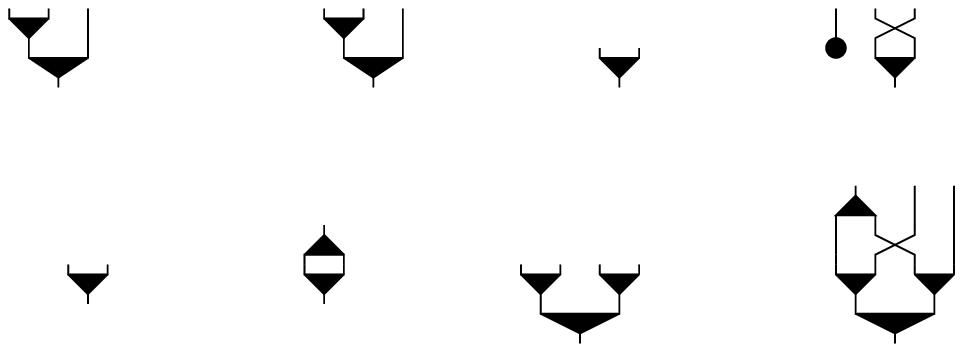}%
\end{picture}%
\setlength{\unitlength}{4144sp}%
\begingroup\makeatletter\ifx\SetFigFont\undefined%
\gdef\SetFigFont#1#2#3#4#5{%
  \reset@font\fontsize{#1}{#2pt}%
  \fontfamily{#3}\fontseries{#4}\fontshape{#5}%
  \selectfont}%
\fi\endgroup%
\begin{picture}(4470,1719)(-137,-883)
\put(3331,524){\makebox(0,0)[b]{\smash{{\SetFigFont{12}{14.4}{\familydefault}{\mddefault}{\updefault}$\Phi$}}}}
\put(  1,704){\makebox(0,0)[b]{\smash{{\SetFigFont{12}{14.4}{\familydefault}{\mddefault}{\updefault}$1$}}}}
\put(271,-466){\makebox(0,0)[b]{\smash{{\SetFigFont{12}{14.4}{\familydefault}{\mddefault}{\updefault}$1$}}}}
\put(451,-466){\makebox(0,0)[b]{\smash{{\SetFigFont{12}{14.4}{\familydefault}{\mddefault}{\updefault}$1$}}}}
\put(181,704){\makebox(0,0)[b]{\smash{{\SetFigFont{12}{14.4}{\familydefault}{\mddefault}{\updefault}$2$}}}}
\put(2881,524){\makebox(0,0)[b]{\smash{{\SetFigFont{12}{14.4}{\familydefault}{\mddefault}{\updefault}$2$}}}}
\put(361,704){\makebox(0,0)[b]{\smash{{\SetFigFont{12}{14.4}{\familydefault}{\mddefault}{\updefault}$3$}}}}
\put(2701,524){\makebox(0,0)[b]{\smash{{\SetFigFont{12}{14.4}{\familydefault}{\mddefault}{\updefault}$3$}}}}
\put(2341,-466){\makebox(0,0)[b]{\smash{{\SetFigFont{12}{14.4}{\familydefault}{\mddefault}{\updefault}$1$}}}}
\put(2701,-466){\makebox(0,0)[b]{\smash{{\SetFigFont{12}{14.4}{\familydefault}{\mddefault}{\updefault}$1$}}}}
\put(2521,-466){\makebox(0,0)[b]{\smash{{\SetFigFont{12}{14.4}{\familydefault}{\mddefault}{\updefault}$2$}}}}
\put(2881,-466){\makebox(0,0)[b]{\smash{{\SetFigFont{12}{14.4}{\familydefault}{\mddefault}{\updefault}$3$}}}}
\put(901,389){\makebox(0,0)[b]{\smash{{\SetFigFont{12}{14.4}{\familydefault}{\mddefault}{\updefault}$\lsur$}}}}
\put(901,-646){\makebox(0,0)[b]{\smash{{\SetFigFont{12}{14.4}{\familydefault}{\mddefault}{\updefault}$\lsur$}}}}
\put(3331,389){\makebox(0,0)[b]{\smash{{\SetFigFont{12}{14.4}{\familydefault}{\mddefault}{\updefault}$\lsur$}}}}
\put(3331,-646){\makebox(0,0)[b]{\smash{{\SetFigFont{12}{14.4}{\familydefault}{\mddefault}{\updefault}$\lsur$}}}}
\put(901,524){\makebox(0,0)[b]{\smash{{\SetFigFont{12}{14.4}{\familydefault}{\mddefault}{\updefault}$\Phi$}}}}
\put(901,-511){\makebox(0,0)[b]{\smash{{\SetFigFont{12}{14.4}{\familydefault}{\mddefault}{\updefault}$\Phi$}}}}
\put(3331,-511){\makebox(0,0)[b]{\smash{{\SetFigFont{12}{14.4}{\familydefault}{\mddefault}{\updefault}$\Phi$}}}}
\end{picture}%
\end{center}

\noindent Then, let us see how to translate the five rules of our term rewriting system derived from the theory of $\zb$-vector spaces. Each rule is pictured in order (associativity, left and right units, commutativity and self-inverse), has been given a name ($A$, $L$, $R$, $C$ and $S$) and has its translation written just below:

\medskip
\begin{center}
\begin{picture}(0,0)%
\includegraphics{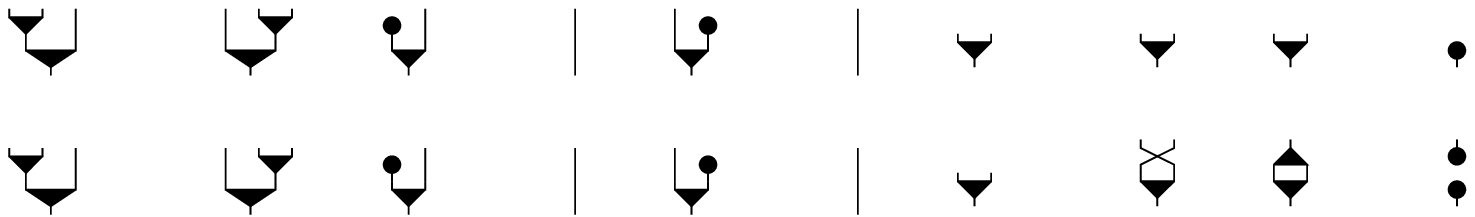}%
\end{picture}%
\setlength{\unitlength}{4144sp}%
\begingroup\makeatletter\ifx\SetFigFont\undefined%
\gdef\SetFigFont#1#2#3#4#5{%
  \reset@font\fontsize{#1}{#2pt}%
  \fontfamily{#3}\fontseries{#4}\fontshape{#5}%
  \selectfont}%
\fi\endgroup%
\begin{picture}(6780,1098)(-135,-438)
\put(6256,324){\makebox(0,0)[b]{\smash{{\SetFigFont{10}{12.0}{\familydefault}{\mddefault}{\updefault}$\red{S}$}}}}
\put(627,324){\makebox(0,0)[b]{\smash{{\SetFigFont{10}{12.0}{\familydefault}{\mddefault}{\updefault}$\red{A}$}}}}
\put(-20,552){\makebox(0,0)[b]{\smash{{\SetFigFont{10}{12.0}{\familydefault}{\mddefault}{\updefault}$1$}}}}
\put(132,552){\makebox(0,0)[b]{\smash{{\SetFigFont{10}{12.0}{\familydefault}{\mddefault}{\updefault}$2$}}}}
\put(284,552){\makebox(0,0)[b]{\smash{{\SetFigFont{10}{12.0}{\familydefault}{\mddefault}{\updefault}$3$}}}}
\put(969,552){\makebox(0,0)[b]{\smash{{\SetFigFont{10}{12.0}{\familydefault}{\mddefault}{\updefault}$1$}}}}
\put(1121,552){\makebox(0,0)[b]{\smash{{\SetFigFont{10}{12.0}{\familydefault}{\mddefault}{\updefault}$2$}}}}
\put(1273,552){\makebox(0,0)[b]{\smash{{\SetFigFont{10}{12.0}{\familydefault}{\mddefault}{\updefault}$3$}}}}
\put(4811,324){\makebox(0,0)[b]{\smash{{\SetFigFont{10}{12.0}{\familydefault}{\mddefault}{\updefault}$\red{C}$}}}}
\put(4317,438){\makebox(0,0)[b]{\smash{{\SetFigFont{10}{12.0}{\familydefault}{\mddefault}{\updefault}$1$}}}}
\put(5306,438){\makebox(0,0)[b]{\smash{{\SetFigFont{10}{12.0}{\familydefault}{\mddefault}{\updefault}$1$}}}}
\put(4469,438){\makebox(0,0)[b]{\smash{{\SetFigFont{10}{12.0}{\familydefault}{\mddefault}{\updefault}$2$}}}}
\put(5153,438){\makebox(0,0)[b]{\smash{{\SetFigFont{10}{12.0}{\familydefault}{\mddefault}{\updefault}$2$}}}}
\put(627,-312){\makebox(0,0)[b]{\smash{{\SetFigFont{10}{12.0}{\familydefault}{\mddefault}{\updefault}$\red{\Phi(A)}$}}}}
\put(4811,-312){\makebox(0,0)[b]{\smash{{\SetFigFont{10}{12.0}{\familydefault}{\mddefault}{\updefault}$\red{\Phi(C)}$}}}}
\put(2224,-312){\makebox(0,0)[b]{\smash{{\SetFigFont{10}{12.0}{\familydefault}{\mddefault}{\updefault}$\red{\Phi(L)}$}}}}
\put(3518,-312){\makebox(0,0)[b]{\smash{{\SetFigFont{10}{12.0}{\familydefault}{\mddefault}{\updefault}$\red{\Phi(R)}$}}}}
\put(6256,-312){\makebox(0,0)[b]{\smash{{\SetFigFont{10}{12.0}{\familydefault}{\mddefault}{\updefault}$\red{\Phi(S)}$}}}}
\put(1882,552){\makebox(0,0)[b]{\smash{{\SetFigFont{10}{12.0}{\familydefault}{\mddefault}{\updefault}$1$}}}}
\put(2567,552){\makebox(0,0)[b]{\smash{{\SetFigFont{10}{12.0}{\familydefault}{\mddefault}{\updefault}$1$}}}}
\put(3518,324){\makebox(0,0)[b]{\smash{{\SetFigFont{10}{12.0}{\familydefault}{\mddefault}{\updefault}$\red{R}$}}}}
\put(3023,552){\makebox(0,0)[b]{\smash{{\SetFigFont{10}{12.0}{\familydefault}{\mddefault}{\updefault}$1$}}}}
\put(3860,552){\makebox(0,0)[b]{\smash{{\SetFigFont{10}{12.0}{\familydefault}{\mddefault}{\updefault}$1$}}}}
\put(5762,438){\makebox(0,0)[b]{\smash{{\SetFigFont{10}{12.0}{\familydefault}{\mddefault}{\updefault}$1$}}}}
\put(5914,438){\makebox(0,0)[b]{\smash{{\SetFigFont{10}{12.0}{\familydefault}{\mddefault}{\updefault}$1$}}}}
\put(2224,324){\makebox(0,0)[b]{\smash{{\SetFigFont{10}{12.0}{\familydefault}{\mddefault}{\updefault}$\red{L}$}}}}
\end{picture}%
\end{center}
\medskip

\noindent Note that several cases may occur. For the first three rules, no resource management operator is added during translation: these three rules are \emph{linear} (or \emph{left-} and \emph{right-linear}). When translated, the commutativity rule has one operator added on its right side and none on its left side: it is a left-linear but not right-linear rule. Finally, the self-inverse rule has one operator added on each of its members during translation: it is neither left- nor right-linear. \\

\noindent Some issues have now been arisen. The first one concerns the rules to be added in order both to describe the behaviour of our local permutation, eraser and duplicator and to ensure the global coherence of these local rules.

The next issue is about the respective computational properties of the starting term rewriting system and of the rewriting system one gets as a result of making the resource management operations explicit. These first two issues are adressed in section \ref{explicit}.

For the moment, we are concerned with a third issue: where does rewriting takes place now? Indeed, starting from a term rewriting system, we have crafted another rewriting system which is not a term one, and for two reasons. The first one is that its signature contains non-algebraic operators, that is operators that do not have exactly one output (the resource management operators have zero or two outputs). The second reason is that variables have been dropped to be replaced by these new operators: this is also a step outside term rewriting. Hence, our new object is not a term rewriting system and section \ref{operads} recalls a notion from [Burroni 1993] used to describe it.

\section{Three-dimensional polygraphs}\label{operads}

\noindent Like equational theories, \emph{$\mathit{3}$-polygraphs} are useful objects in universal algebra, in the sense that they allow one to present algebraic structures by generators and relations. However, they are far more general than equational theories, and this has two consequences: on one hand, they can handle more general objects, like the rewriting system sketched in section \ref{resource}, or the structure of quantum groups; but, on the other hand, their generality comes with an increase in the structural complexity: the development of new tools is mandatory to prove termination, for example. \\

\noindent Polygraphs are genuine categorical objects but we prefer a diagrammatic definition here. For this paper, a $3$-polygraph is made of a signature, that is a set of operators with a finite number of inputs and a finite number of outputs, together with a family of rules: in fact, this is just a special case of $3$-polygraph, one with only one $0$-cell and one $1$-cell. For the complete theory of $n$-polygraphs, the interested reader should check [Burroni 1993].

The operators are once again represented by fixed diagrams of size one, with as many free edges at the top as the operator inputs and as many free edges at the bottom as the operator outputs. For example, some usual diagram shapes are pictured here:
\begin{center}
\begin{picture}(0,0)%
\includegraphics{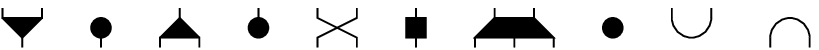}%
\end{picture}%
\setlength{\unitlength}{4144sp}%
\begingroup\makeatletter\ifx\SetFigFont\undefined%
\gdef\SetFigFont#1#2#3#4#5{%
  \reset@font\fontsize{#1}{#2pt}%
  \fontfamily{#3}\fontseries{#4}\fontshape{#5}%
  \selectfont}%
\fi\endgroup%
\begin{picture}(3714,204)(169,467)
\end{picture}%
\end{center}

\noindent Some of them have already been encountered, some of the others are less algebraic: one has zero input and output - it is usefull to describe Petri nets, see [Guiraud 2004] -, one has two inputs and zero output - it is used together with its dual with zero input and two outputs to represent knots and tangles. \\

\noindent Here, the "terms" one considers are all the circuits one can build with all these elementary diagrams: these are the \emph{Penrose diagrams} (or \emph{circuits}) one can build with the size one diagrams representing the operators, such as:
\begin{center}
\begin{picture}(0,0)%
\includegraphics{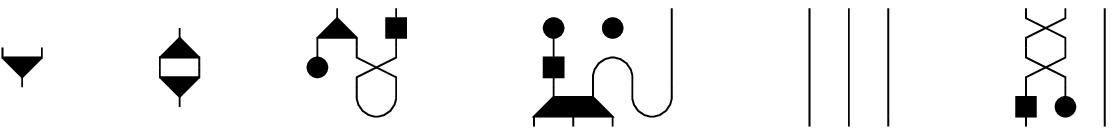}%
\end{picture}%
\setlength{\unitlength}{4144sp}%
\begingroup\makeatletter\ifx\SetFigFont\undefined%
\gdef\SetFigFont#1#2#3#4#5{%
  \reset@font\fontsize{#1}{#2pt}%
  \fontfamily{#3}\fontseries{#4}\fontshape{#5}%
  \selectfont}%
\fi\endgroup%
\begin{picture}(5064,564)(259,-433)
\end{picture}%
\end{center}

\noindent Each of these circuits has a finite number of inputs (on the top) and of outputs (on the bottom) but has no variable. Furthermore, they need not be connected, as the three-inputs and three-outputs wire-only one. \\

\noindent These circuits, which are also called diagrams or arrows, have an algebraic structure. To explain it, let us use the notation $f:m\fl n$ to express that $f$ is a circuit with~$m$ inputs and~$n$ outputs. For any circuit $f$, $s(f)$ is its number of inputs and $t(f)$ its number of outputs. The following constructions and properties are valid for circuits:
\begin{enumerate}
\item[-] Let $f:m\fl n$ and $g:n\fl p$. Then, one can connect each output of $f$ with the corresponding input of $g$, in the same order, to form a new circuit with $m$ inputs and $p$ outputs denoted by $g\circ f$.
\item[-] This composition operation admits local units: a circuit $f:m\fl n$ satisfies $f\circ m=f$ and~$n\circ f=f$, where $p$ is the wire-only circuit with $p$ inputs and $p$ outputs.
\item[-] Let $f:m\fl n$ and $g:p\fl q$. Then, one can put $f$ and $g$ side by side to form a new circuit with $m+p$ inputs and $n+q$ outputs, denoted by $f\tens g$.
\item[-] This product operation admits a bilateral neutral element: the empty circuit $0$ with no input nor output, represented by an empty diagram.
\item[-] Finally, the composition and product are related by the \emph{exchange relations}. They are given by the following equality, that is required to hold for any two circuits $f:m\fl n$ and~$g:p\fl q$:
$$
(t(f)\tens g)\circ (f\tens s(g))=f\tens g=(f\tens t(g))\circ (s(f)\tens g).
$$
\end{enumerate}

\begin{defn}
A family $\Cr$ of circuits endowed with this structure $\tens$ and $\circ$, satisfying the aforegiven unit and exchange relations, is called a \emph{product category}; the subset of circuits with $m$ inputs and $n$ outputs is denoted by $\Cr(m,n)$. When the circuits of $\Cr$ are freely built from a signature $\Sigma$, this object is the \emph{free product category generated by $\Sigma$}, denoted by $\mon{\Sigma}$.
\findef\end{defn}

\begin{rem}
Product categories, or PROs, were defined in [MacLane 1965]. An alternative definition is: a product category is a strict monoidal category whose underlying monoid of objects is $(\Nb,+,0)$, the one of natural numbers with addition and zero. In [Guiraud 2004], such a category was called a \emph{(monochromatic) operad}, for this structure is a common generalization of many universal algebra objects: May's operads, Lawvere's algebraic theories and MacLane's PROs and PROPs.

Product categories are also a special case of \emph{$2$-monoids} or \emph{$2$-categories with only one $0$-cell}. A generalization of this paper results should be possible, since circuit-like diagrams extend to general $2$-cells. For this paper, we stick to MacLane's product categories, but all this terminology will be made clear in subsequent work.
\end{rem}

\noindent A \emph{rewrite rule} on a product category $\Cr$ is a pair $f\fl g$ of \emph{parallel} arrows (they have the same number of inputs and the same number of outputs). Such a rule generates \emph{reductions} on circuits: whenever an arrow $h$ contains $f$, the rule generates a reduction from $h$ to $k$, where $k$ is the same as $h$, except that $f$ has been replaced by $g$. The fact that $f$ and $g$ have the same number of inputs and the same number of outputs ensures that one can connect the unchanged part of the circuit with the changed part, without using implicit operations before.

\begin{defn}
A \emph{$\mathit{3}$-polygraph} is a pair $(\Sigma,R)$ where $\Sigma$ is a signature and $R$ is a family of rewrite rules on $\mon{\Sigma}$.
\end{defn}

\noindent One way to formalize the reduction relation generated by rules on a free product category $\mon{\Sigma}$ is to define \emph{contexts}. We just explain here what they are, avoiding to dig further into the technical aspects, developped in [Guiraud 2004]. Let $\Sigma$ be a signature. Then, a \emph{context} on $\mon{\Sigma}$ is a circuit $c$ with a "hole" inside: this hole has a finite number of inputs and of outputs where on can paste a circuit $f$ with correponding numbers of inputs and outputs; this pasting operation results in a circuit denoted by $c[f]$. Then, a rule $f\fl g$ generates a reduction from each circuit $c[f]$, with $c$ any context, to the circuit $c[g]$. \\

\noindent Finally, given two product categories $\Cr$ and $\Dr$, a \emph{product category functor from $\Cr$ to $\Dr$} is a map which sends each circuit of~$\Cr$ onto a circuit of~$\Dr$ with the same number of inputs and of outputs, and which preserves identities, products and compositions. When $\Cr$ is the free product category $\mon{\Sigma}$, then a classical categorical argument tells us that any product category functor $F:\mon{\Sigma}\fl\Dr$ is entirely and uniquely given by the circuits $F(\phi)$ in $\Dr$, for every operator $\phi$ in $\Sigma$.

\section{From term rewriting to $\mathbf{3}$-polygraphs}\label{explicit}

\noindent This section uses results from [Burroni 1993], presented in a slightly different way, in order to prove a conjecture from [Lafont 2003]: this is theorem \ref{thm0}. This is the result that allows the definition \ref{traduction} of a translation $\Phi$ from any term rewriting system into a $3$-polygraph. Proposition \ref{prop1} and theorem \ref{thm1} give the respective computational properties of the term rewriting system and the $3$-polygraph. \\

\noindent In section \ref{resource}, a $3$-polygraph has been built from the term rewriting system $(\Sigma,R_2)$, which presents the equational theory of $\zb$-vector spaces. Its signature, denoted by $\Sigma^c$, is the one built from $\Sigma$ by addition of the three resource managment operators $\tau$, $\delta$ and $\epsilon$ from $\Delta$. Its family of rules, denoted by $\Phi(R_2)$, consists of the translations $\Phi(A)$, $\Phi(L)$, $\Phi(R)$, $\Phi(C)$ and~$\Phi(S)$ of the five rules from the original term rewriting system. This construction can be generalized to any term rewriting system but is still incomplete for the moment. It lacks two families of rules and this section starts with their description. \\

\noindent Let us fix an algebraic signature $\Sigma$. The set of terms built on the signature $\Sigma$ and on some fixed countable set $V$ of variables is denoted by $T\Sigma$. Let us assume that the set $V$ is endowed with a total order (given by a bijection with $\Nb$), so that the variables can be written $x_1$, $x_2$, $x_3$, etc. For any term~$u$, the notation $\sharp u$ is used for the greatest natural number $i$ such that $x_i$ appears in~$u$. Then, we define $\Tb\Sigma(m,n)$ to be the set of families $(u_1,\dots,u_n)$ of $n$ terms such that $\sharp u_i\leq m$ for every $i$. Note that the set $\Tb\Sigma(m,0)$ has only one element, denoted by~$\ast(m)$. The following operations provide the set $\Tb\Sigma$ with a product category structure:
\begin{enumerate}
\item[-] If $u=(u_1,\dots,u_n)$ is in $\Tb\Sigma(m,n)$ and $v=(v_1,\dots,v_p)$ is in $\Tb\Sigma(n,p)$, then their composite $v\circ u$ is the family $(w_1,\dots,w_p)$ where each $w_i$ is built from $v_i$ by replacing each $x_j$ with $u_j$.
\item[-] The identity of $n$, for any natural number $n$, is the family $(x_1,\dots,x_n)$.
\item[-] The product $u\tens v$ of $u=(u_1,\dots,u_n)$ in $\Tb\Sigma(m,n)$ and of $v=(v_1,\dots,v_q)$ in $\Tb\Sigma(p,q)$ is the family $(w_1,\dots,w_{n+q})$ built that way: if $i$ lies between $1$ and $n$, then $w_i$ is $u_i$; otherwise, $w_{i+n}$ is~$v_i$ where each $x_j$ has been replaced by $x_{j+m}$.
\end{enumerate}

\noindent Furthermore, this product category satisfies some additional properties. The first one is that $\Tb\Sigma$ is a \emph{cartesian} category: seen as a strict monoidal category, the monoidal product $\tens$ is the functorial part of a cartesian product. In our case and informally, this means that every circuit $f:m\fl n$ is entirely and uniquely determined by $n$ circuits $m\fl 1$, in the same way that any function $f:X^m\fl X^n$, where $X$ is a set, is entirely and uniquely determined by $n$ functions $X^m\fl X$: its components. To check that $\Tb\Sigma$ is indeed cartesian, one uses a result from [Burroni 1993], restricted to our setting:

\begin{thm}[Burroni]\label{thm_burroni}
A product category $\Cr$ is cartesian if and only if it contains three arrows:
\begin{center}
\begin{picture}(0,0)%
\includegraphics{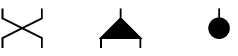}%
\end{picture}%
\setlength{\unitlength}{4144sp}%
\begingroup\makeatletter\ifx\SetFigFont\undefined%
\gdef\SetFigFont#1#2#3#4#5{%
  \reset@font\fontsize{#1}{#2pt}%
  \fontfamily{#3}\fontseries{#4}\fontshape{#5}%
  \selectfont}%
\fi\endgroup%
\begin{picture}(1055,430)(6289,-1469)
\put(7291,-1411){\makebox(0,0)[b]{\smash{{\SetFigFont{12}{14.4}{\familydefault}{\mddefault}{\updefault}{\color[rgb]{0,0,0}$\epsilon$}%
}}}}
\put(6391,-1411){\makebox(0,0)[b]{\smash{{\SetFigFont{12}{14.4}{\familydefault}{\mddefault}{\updefault}{\color[rgb]{0,0,0}$\tau$}%
}}}}
\put(6841,-1411){\makebox(0,0)[b]{\smash{{\SetFigFont{12}{14.4}{\familydefault}{\mddefault}{\updefault}{\color[rgb]{0,0,0}$\delta$}%
}}}}
\end{picture}%
\end{center}

\noindent Such that the two following families of equations hold:
\begin{enumerate}
\item The family $E_{\Delta}$, made of the following seven equations:
\begin{center}
\begin{picture}(0,0)%
\includegraphics{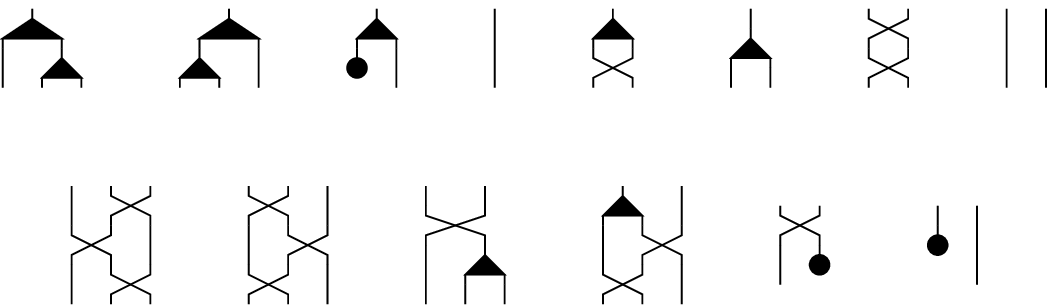}%
\end{picture}%
\setlength{\unitlength}{4144sp}%
\begingroup\makeatletter\ifx\SetFigFont\undefined%
\gdef\SetFigFont#1#2#3#4#5{%
  \reset@font\fontsize{#1}{#2pt}%
  \fontfamily{#3}\fontseries{#4}\fontshape{#5}%
  \selectfont}%
\fi\endgroup%
\begin{picture}(4794,1374)(7639,-5023)
\put(8551,-4786){\makebox(0,0)[b]{\smash{{\SetFigFont{12}{14.4}{\familydefault}{\mddefault}{\updefault}{\color[rgb]{0,0,0}$=$}%
}}}}
\put(8236,-3886){\makebox(0,0)[b]{\smash{{\SetFigFont{12}{14.4}{\familydefault}{\mddefault}{\updefault}{\color[rgb]{0,0,0}$=$}%
}}}}
\put(9676,-3886){\makebox(0,0)[b]{\smash{{\SetFigFont{12}{14.4}{\familydefault}{\mddefault}{\updefault}{\color[rgb]{0,0,0}$=$}%
}}}}
\put(10756,-3886){\makebox(0,0)[b]{\smash{{\SetFigFont{12}{14.4}{\familydefault}{\mddefault}{\updefault}{\color[rgb]{0,0,0}$=$}%
}}}}
\put(12016,-3886){\makebox(0,0)[b]{\smash{{\SetFigFont{12}{14.4}{\familydefault}{\mddefault}{\updefault}{\color[rgb]{0,0,0}$=$}%
}}}}
\put(11656,-4786){\makebox(0,0)[b]{\smash{{\SetFigFont{12}{14.4}{\familydefault}{\mddefault}{\updefault}{\color[rgb]{0,0,0}$=$}%
}}}}
\put(10126,-4786){\makebox(0,0)[b]{\smash{{\SetFigFont{12}{14.4}{\familydefault}{\mddefault}{\updefault}{\color[rgb]{0,0,0}$=$}%
}}}}
\end{picture}%
\end{center}
\item The family $E_{\Sigma}$, made of three equations for each integer $n$ and each arrow $f:n\fl 1$ in $\Cr$:
\begin{center}
\begin{picture}(0,0)%
\includegraphics{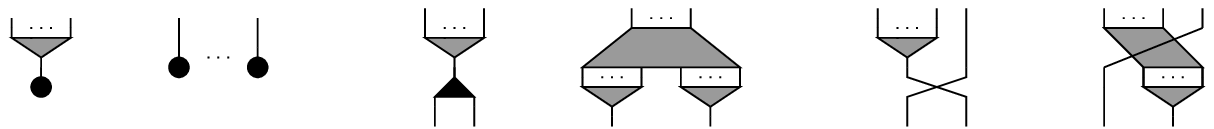}%
\end{picture}%
\setlength{\unitlength}{4144sp}%
\begingroup\makeatletter\ifx\SetFigFont\undefined%
\gdef\SetFigFont#1#2#3#4#5{%
  \reset@font\fontsize{#1}{#2pt}%
  \fontfamily{#3}\fontseries{#4}\fontshape{#5}%
  \selectfont}%
\fi\endgroup%
\begin{picture}(5962,564)(12,107)
\put(5851,164){\makebox(0,0)[b]{\smash{{\SetFigFont{12}{14.4}{\familydefault}{\mddefault}{\updefault}{\color[rgb]{0,0,0}$f$}%
}}}}
\put(766,344){\makebox(0,0)[b]{\smash{{\SetFigFont{12}{14.4}{\familydefault}{\mddefault}{\updefault}{\color[rgb]{0,0,0}$=$}%
}}}}
\put(2656,344){\makebox(0,0)[b]{\smash{{\SetFigFont{12}{14.4}{\familydefault}{\mddefault}{\updefault}{\color[rgb]{0,0,0}$=$}%
}}}}
\put(4951,344){\makebox(0,0)[b]{\smash{{\SetFigFont{12}{14.4}{\familydefault}{\mddefault}{\updefault}{\color[rgb]{0,0,0}$=$}%
}}}}
\put(136,434){\makebox(0,0)[b]{\smash{{\SetFigFont{12}{14.4}{\familydefault}{\mddefault}{\updefault}{\color[rgb]{0,0,0}$f$}%
}}}}
\put(2026,434){\makebox(0,0)[b]{\smash{{\SetFigFont{12}{14.4}{\familydefault}{\mddefault}{\updefault}{\color[rgb]{0,0,0}$f$}%
}}}}
\put(2791,164){\makebox(0,0)[b]{\smash{{\SetFigFont{12}{14.4}{\familydefault}{\mddefault}{\updefault}{\color[rgb]{0,0,0}$f$}%
}}}}
\put(3691,164){\makebox(0,0)[b]{\smash{{\SetFigFont{12}{14.4}{\familydefault}{\mddefault}{\updefault}{\color[rgb]{0,0,0}$f$}%
}}}}
\put(4096,434){\makebox(0,0)[b]{\smash{{\SetFigFont{12}{14.4}{\familydefault}{\mddefault}{\updefault}{\color[rgb]{0,0,0}$f$}%
}}}}
\end{picture}%
\end{center}
\end{enumerate}

\noindent The following recursively defined arrows families $(\delta_n)_{n\in\Nb}$ and $(\tau_{n,1})_{n\in\Nb}$ have been used:
\begin{center}
\begin{picture}(0,0)%
\includegraphics{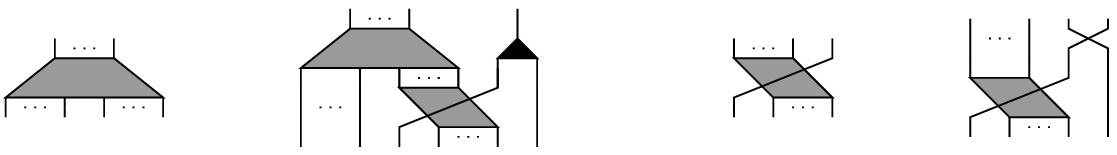}%
\end{picture}%
\setlength{\unitlength}{4144sp}%
\begingroup\makeatletter\ifx\SetFigFont\undefined%
\gdef\SetFigFont#1#2#3#4#5{%
  \reset@font\fontsize{#1}{#2pt}%
  \fontfamily{#3}\fontseries{#4}\fontshape{#5}%
  \selectfont}%
\fi\endgroup%
\begin{picture}(5078,925)(-115,-452)
\put(4006,-16){\makebox(0,0)[b]{\smash{{\SetFigFont{12}{14.4}{\familydefault}{\mddefault}{\updefault}{\color[rgb]{0,0,0}$=$}%
}}}}
\put(271,254){\makebox(0,0)[b]{\smash{{\SetFigFont{8}{9.6}{\familydefault}{\mddefault}{\updefault}{\color[rgb]{0,0,0}$n+1$}%
}}}}
\put( 46,-286){\makebox(0,0)[b]{\smash{{\SetFigFont{8}{9.6}{\familydefault}{\mddefault}{\updefault}{\color[rgb]{0,0,0}$n+1$}%
}}}}
\put(496,-286){\makebox(0,0)[b]{\smash{{\SetFigFont{8}{9.6}{\familydefault}{\mddefault}{\updefault}{\color[rgb]{0,0,0}$n+1$}%
}}}}
\put(3556,-286){\makebox(0,0)[b]{\smash{{\SetFigFont{8}{9.6}{\familydefault}{\mddefault}{\updefault}{\color[rgb]{0,0,0}$n+1$}%
}}}}
\put(3376,254){\makebox(0,0)[b]{\smash{{\SetFigFont{8}{9.6}{\familydefault}{\mddefault}{\updefault}{\color[rgb]{0,0,0}$n+1$}%
}}}}
\put(1621,389){\makebox(0,0)[b]{\smash{{\SetFigFont{8}{9.6}{\familydefault}{\mddefault}{\updefault}{\color[rgb]{0,0,0}$n$}%
}}}}
\put(2026,-421){\makebox(0,0)[b]{\smash{{\SetFigFont{8}{9.6}{\familydefault}{\mddefault}{\updefault}{\color[rgb]{0,0,0}$n$}%
}}}}
\put(1396,-421){\makebox(0,0)[b]{\smash{{\SetFigFont{8}{9.6}{\familydefault}{\mddefault}{\updefault}{\color[rgb]{0,0,0}$n$}%
}}}}
\put(4456,344){\makebox(0,0)[b]{\smash{{\SetFigFont{8}{9.6}{\familydefault}{\mddefault}{\updefault}{\color[rgb]{0,0,0}$n$}%
}}}}
\put(4636,-376){\makebox(0,0)[b]{\smash{{\SetFigFont{8}{9.6}{\familydefault}{\mddefault}{\updefault}{\color[rgb]{0,0,0}$n$}%
}}}}
\put(946,-16){\makebox(0,0)[b]{\smash{{\SetFigFont{12}{14.4}{\familydefault}{\mddefault}{\updefault}{\color[rgb]{0,0,0}$=$}%
}}}}
\end{picture}%
\end{center}

\noindent with the initial values $\delta_0=0$ and $\tau_{0,1}=1$.
\end{thm}

\noindent Note that the following convention is now used in diagrams: generating operators are drawn with black diagrams, while composite arrows are grey. The union of the two families $E_{\Delta}$ and $E_{\Sigma}$ is denoted by~$E_{\Delta\Sigma}$. Theorem \ref{thm_burroni} is not mandatory to get the following proposition but yields an easy proof of it:

\begin{prop}
The product category $\Tb\Sigma$ is cartesian.
\end{prop}

\begin{dem}
Let us start with the definition of the three arrows from theorem \ref{thm_burroni}: the arrow $\tau$ is the pair $(x_2,x_1)$ of terms; the arrow $\delta$ is $(x_1,x_1)$; finally, the arrow $\epsilon$ is the empty family $\ast(1)$. Computations to check the equations of theorem \ref{thm_burroni} are straightforward.
\findem\end{dem}

\noindent The next step consists in the proof that $\Tb\Sigma$ is the \emph{free} cartesian category generated by the algebraic signature $\Sigma$. In order to prove this fact, one starts with another use of theorem \ref{thm_burroni}:

\begin{cor}[of theorem \ref{thm_burroni}]\label{cor_burroni}
For every algebraic signature $\Sigma$, the category $\mon{\Sigma^c}/E_{\Delta\Sigma}$ is the free cartesian category generated by $\Sigma$.
\end{cor}

\noindent Hence, in order to prove that $\Tb\Sigma$ is another version of the free cartesian category generated by $\Sigma$, it is sufficient to prove that there exists an isomorphism $\hat{\Phi}:\Tb\Sigma\fl\mon{\Sigma^c}/E_{\Delta\Sigma}$.

The signature $\Sigma$ is contained in $\Tb\Sigma$: one defines an inclusion $i$ which sends each $\phi:n\fl 1$ from $\Sigma$ onto the term $\phi(x_1,\dots,x_n)$. Hence, corollary \ref{cor_burroni} extends $i$ into a cartesian functor $F$ from $\mon{\Sigma^c}/E_{\Delta\Sigma}$ to $\Tb\Sigma$: this functor sends each $\phi$ from $\Sigma$ onto $i(\phi)$ and $\tau$, $\delta$ and~$\epsilon$ respectively onto $(x_2,x_1)$, $(x_1,x_1)$ and~$\ast(1)$.

Conversely, let us consider an arrow $f=(u_1,\dots,u_n)$ in $\Tb\Sigma(m,n)$. Each term $u_i$ can be written $u_i=f_i(y^i_1,\dots,y^i_{k_i})$, with $k_i$ an integer, $f_i$ an arrow in $\mon{\Sigma}(k_i,1)$ and each $y^i_j$ a variable from $\ens{x_1}{x_m}$. Furthermore, this decomposition of terms is unique. Thus, the arrow $f$ uniquely decomposes into:
$$
f=(f_1\tens\dots\tens f_n)\circ(y^1_1,\dots,y^n_{k_n}).
$$

\noindent There remains to prove that every family $(y_1,\dots,y_k)$ of variables in $\ens{x_1}{x_m}$ can be uniquely written (\emph{modulo} $E_{\Delta}$) with the three arrows $i(\tau)$, $i(\delta)$ and $i(\epsilon)$. This can be done in two steps.

Let us define the sub-product category $\Vb$ of $\Tb\Sigma$ by restricting ourselves to families of variables: this is $\Tb\emptyset$, where $\emptyset$ denotes the signature with no operator. One also defines the cartesian category $\Fb^o$ of \emph{finite sets} with: the arrows of $\Fb^o(m,n)$ are in bijective correspondance with the functions from the finite set $[n]=\ens{1}{n}$ to $[m]$. Then:

\begin{lem}
The cartesian categories $\Vb$ and $\Fb^o$ are isomorphic.
\end{lem}

\begin{dem}
Let $(y_1,\dots,y_n)$ be a family of variables taken in $\ens{x_1}{x_m}$. Then, there exists an unique function $f^*$ from $[n]$ to $[m]$ such that $y_i=x_{f^*(i)}$ for each $i$. Let us fix $\theta(y_1,\dots,y_n)$ as the arrow $f$ in~$\Fb^o$ that corresponds to $f^*$. Conversely, if $f$ is an arrow in $\Fb^o(m,n)$: let us denote by $f^*$ the corresponding function from $[n]$ to $[m]$. Then one defines $\omega(f)=(x_{f^*(1)},\dots,x_{f^*(n)})$. There remains to check that $\theta$ and $\omega$ are cartesian functors which are inverse one another, which is straightforward.
\findem\end{dem}

\noindent The second step uses another result from [Burroni 1993]:

\begin{thm}[Burroni]
The cartesian categories $\Fb^o$ and $\mon{\Delta}/E_{\Delta}$ are isomorphic.
\end{thm}

\noindent Hence, the cartesian categories $\Vb$ and $\mon{\Delta}/E_{\Delta}$ are isomorphic. Consequently, each family $(y_1,\dots,y_k)$ of variables taken in $\ens{x_1}{x_n}$ corresponds to a unique arrow in $\mon{\Delta}/E_{\Delta}$. Furthermore, each arrow $f$ in $\Tb\Sigma(m,n)$ admits a unique decomposition $f=f_{\Sigma}\circ f_{\Delta}$ with $f_{\Sigma}$ in $\mon{\Sigma}$ and $f_{\Delta}$ in $\Vb$.

Finally, one gets that the cartesian functor $F$ from $\mon{\Sigma^c}/E_{\Delta\Sigma}$ to $\Tb\Sigma$ is an isomorphism. However, we want an map from $\Tb\Sigma$ to $\mon{\Sigma^c}$: let us find a convergent $3$-polygraph $(\Sigma^c,R_{\Delta\Sigma})$ such that $\mon{\Sigma^c}/R_{\Delta\Sigma}$ is isomorphic to $\mon{\Sigma^c}/E_{\Delta\Sigma}$ and use the unique normal form property.

A conjecture from [Lafont 2003] is proved:

\begin{thm}\label{thm0}
For any algebraic signature $\Sigma$, the $3$-polygraph $(\Sigma^c,R_{\Delta\Sigma})$ is convergent and $\mon{\Sigma^c}/R_{\Delta\Sigma}$ is isomorphic to the free cartesian category $\mon{\Sigma^c}/E_{\Delta\Sigma}$ generated by $\Sigma$, where the family of rules $R_{\Delta\Sigma}$ is made of the following two subfamilies:
\begin{enumerate}
\item The family $R_{\Delta}$:
\begin{center}
\begin{picture}(0,0)%
\includegraphics{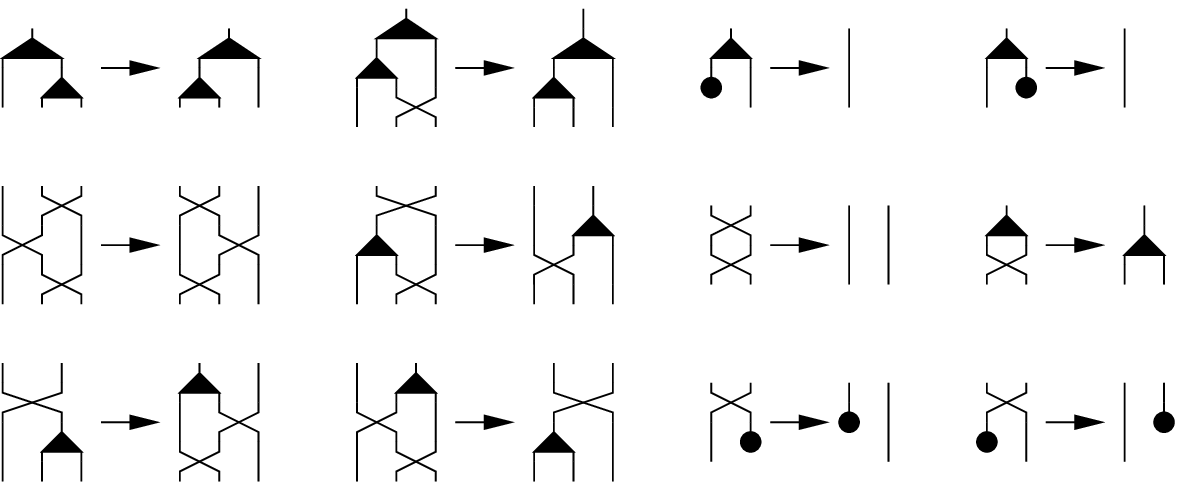}%
\end{picture}%
\setlength{\unitlength}{4144sp}%
\begingroup\makeatletter\ifx\SetFigFont\undefined%
\gdef\SetFigFont#1#2#3#4#5{%
  \reset@font\fontsize{#1}{#2pt}%
  \fontfamily{#3}\fontseries{#4}\fontshape{#5}%
  \selectfont}%
\fi\endgroup%
\begin{picture}(5375,2184)(79,-1423)
\end{picture}%
\end{center}

\item The family $R_{\Sigma}$ given, for each integer $n$ and each operator $\phi$ in $\Sigma(n,1)$, by:
\begin{center}
\begin{picture}(0,0)%
\includegraphics{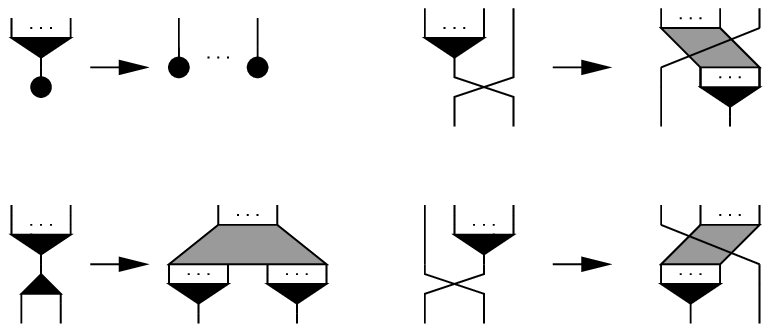}%
\end{picture}%
\setlength{\unitlength}{4144sp}%
\begingroup\makeatletter\ifx\SetFigFont\undefined%
\gdef\SetFigFont#1#2#3#4#5{%
  \reset@font\fontsize{#1}{#2pt}%
  \fontfamily{#3}\fontseries{#4}\fontshape{#5}%
  \selectfont}%
\fi\endgroup%
\begin{picture}(3690,1510)(136,-839)
\put(3151,-781){\makebox(0,0)[b]{\smash{{\SetFigFont{12}{14.4}{\familydefault}{\mddefault}{\updefault}{\color[rgb]{0,0,0}$\phi$}%
}}}}
\put(136,434){\makebox(0,0)[b]{\smash{{\SetFigFont{12}{14.4}{\familydefault}{\mddefault}{\updefault}{\color[rgb]{0,0,0}$\phi$}%
}}}}
\put(136,-466){\makebox(0,0)[b]{\smash{{\SetFigFont{12}{14.4}{\familydefault}{\mddefault}{\updefault}{\color[rgb]{0,0,0}$\phi$}%
}}}}
\put(946,-781){\makebox(0,0)[b]{\smash{{\SetFigFont{12}{14.4}{\familydefault}{\mddefault}{\updefault}{\color[rgb]{0,0,0}$\phi$}%
}}}}
\put(1711,-781){\makebox(0,0)[b]{\smash{{\SetFigFont{12}{14.4}{\familydefault}{\mddefault}{\updefault}{\color[rgb]{0,0,0}$\phi$}%
}}}}
\put(2026,434){\makebox(0,0)[b]{\smash{{\SetFigFont{12}{14.4}{\familydefault}{\mddefault}{\updefault}{\color[rgb]{0,0,0}$\phi$}%
}}}}
\put(2656,-421){\makebox(0,0)[b]{\smash{{\SetFigFont{12}{14.4}{\familydefault}{\mddefault}{\updefault}{\color[rgb]{0,0,0}$\phi$}%
}}}}
\put(3826,209){\makebox(0,0)[b]{\smash{{\SetFigFont{12}{14.4}{\familydefault}{\mddefault}{\updefault}{\color[rgb]{0,0,0}$\phi$}%
}}}}
\end{picture}%
\end{center}
\end{enumerate}
\end{thm}

\begin{rem}
Three families of verifications need to be done. The first one consists in checking that the new rules are derivable from $E_{\Delta\Sigma}$, which is straightforward.

The second one is much more complicated: one needs to check that the $3$-polygraph terminates. However, the structural complexity of polygraphs requires new techniques since the usual ones used in rewriting do not work. One way to craft reduction orders for $3$-polygraphs is made explicit in section~\ref{technique} and used in section \ref{app1} in order to prove the termination of $(\Sigma^c,R_{\Delta\Sigma})$.

Finally, one needs to check that this $3$-polygraph is confluent. Here, this is equivalent to computing all of its critical pairs and check that each one is confluent. Once again, the structural complexity of polygraphs generates problems unknown with other kinds of rewriting theories. For example, a finite $3$-polygraph can produce an infinite number of critical pairs; this is the case here. However, among these critical pairs, some have properties that allow us to finally have only a finite number of computations to do. Critical pairs of $3$-polygraphs need to be further studied and classified according to properties of this kind; this will be addressed in subsequent work.

The present case is discussed in section \ref{app1} and fully studied in [Guiraud 2004].
\end{rem}

\noindent From theorem \ref{thm0}, one concludes the existence of a map $\Phi$ from $\Tb\Sigma$ to $\mon{\Sigma^c}$. Indeed, if $f$ is an arrow in the cartesian category $\Tb\Sigma$, then $\Phi(f)$ will be the $R_{\Delta\Sigma}$-normal form of any representant in $\mon{\Sigma^c}$ of the arrow $F(f)$ in the product category $\mon{\Sigma^c}/E_{\Delta\Sigma}$. This map $\Phi$, which could not be proved to exist until theorem \ref{thm0}, allows the formal definition of the translation of terms into circuits.

\begin{defn}\label{traduction}
For every term $u$ in $T\Sigma$ and for every integer $n\geq\sharp u$, the term $u$ can be seen as an arrow $u_n$ in $\Tb\Sigma(n,1)$. One denotes by $\Phi^n(u)$ the arrow $\Phi(u_n)$ of $\mon{\Sigma^c}$ and by $\Phi(u)$ the particular case~$\Phi^{\sharp u}(u)$. If $\alpha=(u,v)$ is a rewrite rule on $T\Sigma$, the notation $\Phi(\alpha)$ is used for the rewrite rule~$(\Phi(u),\Phi^{\sharp u}(v))$ on $\mon{\Sigma^c}$.
\end{defn}

\noindent As an immediate consequence of the definition, one gets:

\begin{lem}
For any algebraic signature $\Sigma$, any term $u$ in $T\Sigma$ and any integer $n\geq\sharp u$, the arrow $\Phi^n(u)$ is a normal form for the resource management rules $R_{\Delta\Sigma}$.
\end{lem}

\noindent The rest of this section is devoted to the comparison of a term rewriting system $(\Sigma,R)$ with the $3$-polygraph $(\Sigma^c,R^c)$, where $R^c$ is the union of the family $R_{\Delta\Sigma}$ of resource management rules and of the family $\Phi(R)$ made of the translations by $\Phi$ of the rules $R$.

\begin{rem}
Before stating the result, let us qualify by \emph{uniformized} a rule $(u,v)$ on $T\Sigma$ such that $u=f(y_1,\dots,y_k)$ with $f$ an arrow in $\mon{\Sigma}$ and $(y_1,\dots,y_k)$ a family of variables with the following property: $y_1$ is $x_1$; then, for each $i$ in $\ens{1}{k-1}$, the variable $y_{i+1}$ is either in $\ens{y_1}{y_i}$, or~$y_{i+1}$ is~$x_{p+1}$ if $\ens{y_1}{y_i}=\ens{x_1}{x_p}$.

Note that any rule on $T\Sigma$ can be replaced by a uniquely defined uniformized rule that generates the same reduction relation. Furthermore, if a left-linear rule is replaced by its uniformized rule, this one is also left-linear.

Hence, for what follows, (left-linear) term rewriting systems can always be considered uniformized: if they are not, they are replaced by their uniformized equivalent version, with no consequence on the results.

This choice simplifies the translations: a rule $(u,v)$ that is both left-linear and uniformized satisfies $u=f(x_1,\dots,x_{\sharp u})$, with $f$ an arrow in $\mon{\Sigma}$, uniquely defined; hence, the translation by $\Phi$ of such a $u$ is~$f$ and thus is an arrow of $\mon{\Sigma}$.
\end{rem}

\begin{prop}\label{prop1}
If $(\Sigma,R)$ is a term rewriting system, then:
\begin{enumerate}
\item[1.] If the term rewriting system $(\Sigma,R)$ terminates, so does the $3$-polygraph $(\Sigma^c,R^c)$.
\item[2.] The translation $\Phi$ preserves the reduction steps generated by any left-linear rule $\alpha$, that is: for any pair $(u,v)$ of terms such that $u\red{\alpha}\:v$ and any integer $n\geq\sharp u$, there exists an arrow~$f$ in $\mon{\Sigma^c}$ such that
$$
\Phi^n(u)\red{\Phi(\alpha)}\:f\mred{R_{\Delta\Sigma}}\Phi^n(v).
$$
\end{enumerate}
\end{prop}

\begin{dem}
Point 1 uses the technique to be introduced in section \ref{technique}. Its proof is thus postponed until section~\ref{app1}. Point 2 requires lengthy and cumbersome though intuitively simple computations that can be found in [Guiraud 2004].
\findem\end{dem}

\begin{thm}\label{thm1}
A left-linear term rewriting system $(\Sigma,R)$ terminates (resp. is confluent) if and only if its associated $3$-polygraph $(\Sigma^c,R^c)$ terminates (resp. is confluent).
\end{thm}

\begin{dem}
Let us assume that the $3$-polygraph $(\Sigma^c,R^c)$ terminates while the term rewriting system $(\Sigma,R)$ does not. Consequently, there exists some sequence $(u_n)_{n\in\Nb}$ of terms in $T\Sigma$ such that $u_n\red{R}\:u_{n+1}$ for every $n$. From \ref{prop1}, since every rule in $R$ is left-linear, one concludes that, for every $k\geq\sharp u_0$:
$$
\Phi^k(u_0)\pmred{R^c}\:\Phi^k(u_1)\pmred{R^c}\:\cdots\pmred{R^c}\:\Phi^k(u_n)\pmred{R^c}\:\Phi^k(u_{n+1})\pmred{R^c}\:\cdots
$$

\noindent where the notation $\pmred{R^c}$ stands for a \emph{non-empty} $R^c$-reduction path. Such an infinite reduction path existence is denied by the termination of the $3$-polygraph $(\Sigma^c,R^c)$, thus giving this property for $(\Sigma,R)$. The converse, which is true even if the term rewriting system is not left-linear, is still postponed to section \ref{app1}. \\

\noindent Now, let us assume that the term rewriting system $(\Sigma,R)$ is confluent. Let us consider a branching $(f,g,h)$ of $(\Sigma^c,R^c)$: the arrows $f$, $g$ and $h$ have the same finite number of inputs, say $m$, and the same finite number of outputs, say $n$, and satisfy $f\mred{R^c}\:g$ and $f\mred{R^c}\:h$. Let us denote by $\pi$ the canonical projection of $\mon{\Sigma^c}$ onto $\Tb\Sigma$. Then, $\pi$ sends each of $f$, $g$ and $h$ on families $(f_1,\dots,f_n)$, $(g_1,\dots,g_n)$ and $(h_1,\dots,h_n)$ of terms such that each one has variables in $\ens{x_1}{x_m}$. Moreover, for each $i$, one gets that the triple $(f_i,g_i,h_i)$ is a branching of $(\Sigma,R)$. From confluence of this rewriting system, one concludes the existence of a arrow $k_i$ that closes this branching. Let us define $k$ as the translation, by~$\Phi^m$, in $\mon{\Sigma^c}(m,n)$, of the family $(k_1,\dots,k_n)$ of terms. Since $(\Sigma,R)$ is left-linear, proposition \ref{prop1} ensures that this arrow $k$ closes the branching $(f,g,h)$.

Conversely, let us assume that the $3$-polygraph $(\Sigma^c,R^c)$ is confluent. Let us consider a branching $(u,v,w)$ in $(\Sigma,R)$; since this rewriting system is left-linear, this branching translates to a branching $(\Phi^n(u),\Phi^n(v),\Phi^n(w))$ in $(\Sigma^c,R^c)$ for any $n\geq\sharp u$. Since the $3$-polygraph is confluent, there exists some arrow $f$ in $\mon{\Sigma^c}(n,1)$ closing this branching. The projection $\pi(f)$ is an arrow in $\Tb\Sigma(n,1)$ and thus corresponds to a term that closes the initial branching $(u,v,w)$.
\findem\end{dem}

\noindent Before considering what this result allows (or rather does not allow) us to conclude about our term rewriting system $(\Sigma,R_2)$ presenting the theory of $\zb$-vector spaces, there remains some termination results to prove in the elapsed section. However, the intrinsic complexity of the polygraph structure prevents the use of classical techniques; rather, the incoming section presents an adaptation to the particular case of $3$-polygraph we consider of classical interpretation techniques used to craft termination orders for terms.

\section{Termination orders for $\mathbf{3}$-polygraphs}\label{technique}

\noindent In rewriting, one of the most used technique to prove termination is the following one: build a \emph{reduction order}, which is a terminating strict order that is compatible with the term structure; then prove that this order contains the rules. Hence any reduction path in the corresponding rewriting system yields a strictly decreasing family for the reduction order: the fact that such families cannot be infinite ensures that there cannot exist any infinite reduction path or, equivalently, that the rewriting system terminates.

In term rewriting, one easy way to build reduction orders is by means of an \emph{interpretation}. The simpliest ones are: each term $u$ such that $\sharp u=n$ is sent to a function $u_*$ from $\Nb^n$ to $\Nb$ (or any set equipped with a terminating strict order). Then, one says that $u>v$ if each $n$-uple of integers is sent to a strictly greater integer by $u_*$ than by $v_*$. One easy way to compute $u_*$ for each term $u$ is to fix $\phi_*$ for each operator $\phi$ in the considered signature and to extend these values functorially. If one can prove that each~$\phi_*$ is a strictly monotone map and that $f>g$ for each rule $f\fl g$, then $u>v$ whenever there is a reduction from $u$ to $v$. Since the order on $\Nb$ is terminating, so is the order on functions: hence, the considered term rewriting system terminates.

However, in the case of $3$-polygraphs, this classical interpretation technique does not yield reduction orders in general. Indeed, it is not always possible to send each operator $\phi$ of the signature onto a strictly monotone map: for example, the erasure operator $\epsilon$ will be sent to an function from $\Nb$ to $\Nb^0$, that is to a single-element set: this function is unique and monotone, but not strictly. Consequently: even if a rule $f\red\:g$ satisfies $f_*>g_*$, then $(\epsilon^n\circ f)_*=(\epsilon^n\circ g)_*$, with $n$ the number of outputs of both $f$ and~$g$.

One could also consider contravariant interpretations: hence, $\epsilon$ would be sent to a constant natural number~$\epsilon^*$. But, in the most interesting $3$-polygraphs, such as the ones we are concerned with, there is a constant operator $\eta$ which cannot be contravariantly sent to a strictly monotone map. The interpretation technique must be adaptated to the polygraph structure in order to yield termination orders. \\

\noindent Here we are in front of a choice between two possible directions: the first one consists in interpreting arrows into functions between objects equipped with a monoidal product, rather than a cartesian one, such as vector spaces. But, when examined, this has led to horrendous computations that did not produce any reduction order. Nonetheless, this trail is not to be forgotten and shall be reexamined when there is a computational tool, adaptated to polygraphs.

The other path consists in using classical interpretations, both covariant and contravariant, as tools to build a third interpretation: this one will give the desired reduction orders. Let us present images that describe the intuition beneath the formalism. Each arrow in the considered product category is seen as an electrical circuit whose elementary components are the operators it is built from, such as suggested by the diagrammatic representation used. Then, a heat production value is associated to each circuit: each of its inputs and outputs receives a current with a fixed intensity; hence there are two types of currents: some are descending (they come from the inputs and propagate downwards to the outputs) and some are ascending (they propagate upwards, from the outputs to the inputs).

The heat produced by a fixed circuit is calculated this way: an operator is arbitrarily chosen. Then, currents are propagated through the other operators to the chosen one. This requires that choices have been made for each operator: for each one, one must be able to compute the intensities of descending currents transmitted when he knows the intensities of incoming descending current, and similarly with ascending currents. When one knows the intensities of each current coming into the chosen operator, one computes the heat it produces, according to values fixed in advance. Then, one repeats the same procedure for each operator, and sums the results to get the heat produced by the considered circuit, for the chosen current intensities.

Two circuits with the same number of inputs and the same number of outputs are compared this way: if, for each family of (ascending and descending) current intensities, one produces more heat than the other one, then the first one is said to be greater. The goal of this section is twofolds: firstly, to formalize the objects required to compute such an order; secondly, to obtain sufficient conditions for this order to be a reduction order. \\

\noindent Let us describe the required materials. The first one is the object where the interpretations take their values: this will be a product category equipped with a strict order. In order to build it, one considers (non-empty) ordered sets $X$ and $Y$ to express the current intensities, one for descending currents, one for ascending currents (for one of the applications to be described, two different sets of values are needed). Then, a commutative monoid $M$ will contain the possible values of heats; moreover, it is supposed to be equipped with an order such that the addition is strictly monotone in both variables.

From the data $X$, $Y$ and $M$, one builds a somewhat weird product category $\Or(X,Y,M)$ this way: an arrow from $m$ to $n$ in $\Or(X,Y,M)$ is a triple $f=(f_*,f^*,[f])$ consisting of three \emph{monotone} functions
$$
f_*:X^m\fl X^n, \quad f^*:Y^n\fl Y^m, \quad [f]:X^m\times Y^n\fl M.
$$

\noindent The identity of $n$ is the triple $n=(X^n,Y^n,0)$ made from the identities of $X^n$ and $Y^n$ and the constant zero-function from $X^n\times Y^n$ to $M$. Two arrows $f:m\fl n$ and $g:n\fl p$ compose this way: $(g\circ f)_*$ and $(g\circ f)^*$ are respectively the composites $g_*\circ f_*$ and $f^*\circ g^*$; for elements $\vec{x}$ in $X^m$ and $\vec{y}$ in $Y^p$, the function $[g\circ f]$ is given by:
$$
[g\circ f](\vec{x},\vec{y})=[f](\vec{x},g^*(\vec{y}))+[g](f_*(\vec{x}),\vec{y}).
$$

\noindent If $f:m\fl n$ and $g:p\fl q$ are two arrows in $\Or(X,Y,M)$, then their product is given by: $(f\tens g)_*$ and $(f\tens g)^*$ are respectively $f_*\tens g_*$ and $f^*\tens g^*$; if $\vec{x}$, $\vec{y}$, $\vec{x}'$ and $\vec{y}\,'$ are respectively elements of $X^m$, $Y^n$, $X^p$ and $Y^q$, then $[f\tens g]$ is given by:
$$
[f\tens g](\vec{x},\vec{x}',\vec{y},\vec{y}\,')=[f](\vec{x},\vec{y})+[g](\vec{x}',\vec{y}\,').
$$

\noindent Then one checks that these operations return monotone functions and that they satisfy the required equations, in order to get:

\begin{lem}
The aforedefined object $\Or(X,Y,M)$ is a product category.
\end{lem}

\noindent On top of this product category structure, a strict order relation $\succ$ is defined on parallel arrows of $\Or(X,Y,M)$. If $f$ and $g$ are two arrows from $m$ to $n$, then $f\succ g$ if, for any $\vec{x}$ in $X^m$ and $\vec{y}$ in $Y^n$, the following three inequalities hold:
$$
f_*(\vec{x})\geq g_*(\vec{x}), \quad f^*(\vec{y})\geq g^*(\vec{y}), \quad [f](\vec{x},\vec{y})>[g](\vec{x},\vec{y}).
$$

\noindent Now, let us consider a signature $\Sigma$. Let us asume that each operator $\phi:m\fl n$ in $\Sigma$ is associated with an arrow $(\phi_*,\phi^*,[\phi]):m\fl n$ in $\Or(X,Y,M)$: this is the interpretation. For any $\phi$, the monotone functions $\phi_*$, $\phi^*$ and $[\phi]$ respectively express how the operator transmits descending and ascending currents and how much heat it produces, according to the current intensities it receives.

Since $\mon{\Sigma}$ is the free product category generated by the signature $\Sigma$, the map sending each $\phi$ in $\Sigma$ to the triple $(\phi_*,\phi^*,[\phi])$ uniquely extends to a product category functor $F$ from $\mon{\Sigma}$ to $\Or(X,Y,M)$. This means that one can compute $f_*$, $f^*$ and $[f]$ for any circuit~$f$ in~$\mon{\Sigma}$, from the values $\phi_*$, $\phi^*$ and $[\phi]$ given for each operator $\phi$ in $\Sigma$ and using the formulas for composition and product in $\Or(X,Y,M)$.

The last step consists in using $F$ to get the order $\succ$ back from $\Or(X,Y,M)$ on $\mon{\Sigma}$: for any two parallel arrows $f$ and $g$ in $\mon{\Sigma}$, then $f\succ g$ is $F(f)\succ F(g)$.

\begin{thm}\label{thm2}
With the aforegiven notations and if the strict part of the order on $M$ is terminating, then the strict order~$\succ$ constructed on $\mon{\Sigma}$ is a reduction order.
\end{thm}

\begin{dem}
One must check that the binary relation $\succ$ built on $\mon{\Sigma}$ is antireflexive, transitive, terminating and compatible with the product category structure. Let us assume that $f$ is an arrow in $\mon{\Sigma}(m,n)$ such that $f\succ f$; let us fix any elements $\vec{x}$ and $\vec{y}$ respectively in the non-empty sets $X^m$ and $Y^n$; then, by definition of $\succ$, one gets the following strict inequality in $M$:
$$
[F(f)](\vec{x},\vec{y})>[F(f)](\vec{x},\vec{y}).
$$

\noindent However, this inequality cannot hold in $M$ since $>$ is the strict part of an order relation. The termination is proved by a similar argument: any infinite and strictly decreasing sequence in $\mon{\Sigma}$ yields, through the non-emptyness of $X$ and $Y$, at least one infinite strictly decreasing sequence in $M$, which existence is denied by the assumed termination of the strict part of its order. The transitivity comes from the ones of the orders on $X$, $Y$ and $M$. Finally, compatibility with the product category structure is checked through computations which use the monotone quality of each $f_*$, $f^*$ and $[f]$ in $\Or(X,Y,M)$, together with the facts that $M$ is a commutative monoid and $F$ is an product category functor.
\findem\end{dem}

\noindent For concrete applications, presented in the next two sections, the following corollary will be used instead of theorem \ref{thm2}:

\begin{cor}\label{cor2}
Let us consider a $3$-polygraph $(\Sigma,R)$. Let us assume that there exist:
\begin{enumerate}
\item Two non-empty ordered sets $X$ and $Y$.
\item A commutative monoid $M$ equipped with an order such that its strict part is \emph{terminating} and such that the sum is strictly monotone in both variables.
\item For each operator $\phi$ in $\Sigma(m,n)$, three \emph{monotone} functions:
$$
\phi_*:X^m\fl X^n, \quad \phi^*:Y^n\fl Y^m, \quad [\phi]:X^m\times Y^n\fl M.
$$
\end{enumerate}

\noindent If the strict order $\succ$ on arrows of $\mon{\Sigma}$ built from these data, in the aforegiven manner, satisfies $f\succ g$ for every rule $f\fl g$ in $R$, then the $3$-polygraph $(\Sigma,R)$ terminates.
\end{cor}

\section{Application 1: explicit resource management polygraphs}\label{app1}

\noindent This section is devoted to the remaining unproved results from section \ref{explicit}. Let us fix a term rewriting system $(\Sigma,R)$ for the whole section.

\subsection{Convergence of the $\mathbf{3}$-polygraph of explicit resource management}

\noindent The first result to prove is theorem \ref{thm0}: the $3$-polygraph $(\Sigma^c,R_{\Delta\Sigma})$ is convergent, where we recall from section~\ref{explicit} that $\Sigma^c$ is the signature made of the algebraic signature $\Sigma$ and the resource management signature $\Delta$, while $R_{\Delta\Sigma}$ is the family of resource management rules.

The proof is divided in three steps: the first one consists in proving its termination; then, we recall from [Guiraud 2004] that this $3$-polygraph is locally confluent; finally, Newman's lemma is applied to get its convergence. Let us start with termination: we use the technique developped in section \ref{technique}. However, the considered polygraph is rather complex and needs two applications of the technique. For the rest of this paragraph, let us fix some notations. Let us denote by $\alpha$ the following rule:
\begin{center}
\begin{picture}(0,0)%
\includegraphics{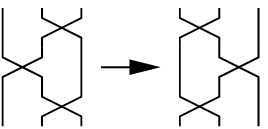}%
\end{picture}%
\setlength{\unitlength}{4144sp}%
\begingroup\makeatletter\ifx\SetFigFont\undefined%
\gdef\SetFigFont#1#2#3#4#5{%
  \reset@font\fontsize{#1}{#2pt}%
  \fontfamily{#3}\fontseries{#4}\fontshape{#5}%
  \selectfont}%
\fi\endgroup%
\begin{picture}(1194,564)(79,-613)
\end{picture}%
\end{center}

\noindent We denote by $\Nb^*$ the set of non-zero natural numbers with its natural order relation. The commutative monoid freely generated by $\Nb^*$ is denoted by $[\Nb^*]$ and is considered equipped by the \emph{multiset order} generated by the usual order relation on natural numbers. The elements of $[\Nb^*]$ are all the finite formal sums of non-zero natural numbers; a natural number $n$, seen as a generator of $[\Nb^*]$, is denoted by $\ul{n}$.

\noindent The multiset order is defined in two steps: for the first one, one says that any sum $a=\sum_i k_i.\ul{n_i}$ satisfies the inequality $\ul{n}>a$ if $n>n_i$ for each $i$; then, the multiset order is taken as the reflexive and structure-compatible closure of this relation.

This implies that the addition is strictly monotone in both variables; furthermore, since the strict order $>$ on $\Nb^*$ terminates, so does the strict part of the multiset order. Here is an example of some strict inequalities that hold in $[\Nb^*]$:
$$
0 \: < \: 127.\ul{1} \:<\: \ul{2} \:<\: 4.\ul{1}+2.\ul{3} \:<\: \ul{4}\;.
$$

\begin{lem*}\label{lem1}
The $3$-polygraph $(\Sigma^c,R_{\Delta\Sigma})$ terminates if and only if the $3$-polygraph $(\Sigma^c,\{\alpha\})$ terminates.
\end{lem*}

\begin{dem}
Let us consider the product category $\Or(\Nb^*,\Nb^*,[\Nb^*])$ together with the termination order $\succ$ as defined in section \ref{technique}. Let us denote by $F$ the product category functor from $\mon{\Sigma^c}$ into $\Or(\Nb^*,\Nb^*,[\Nb^*])$ defined by the following values on the operators of $\Sigma^c$:
\begin{center}
\begin{picture}(0,0)%
\includegraphics{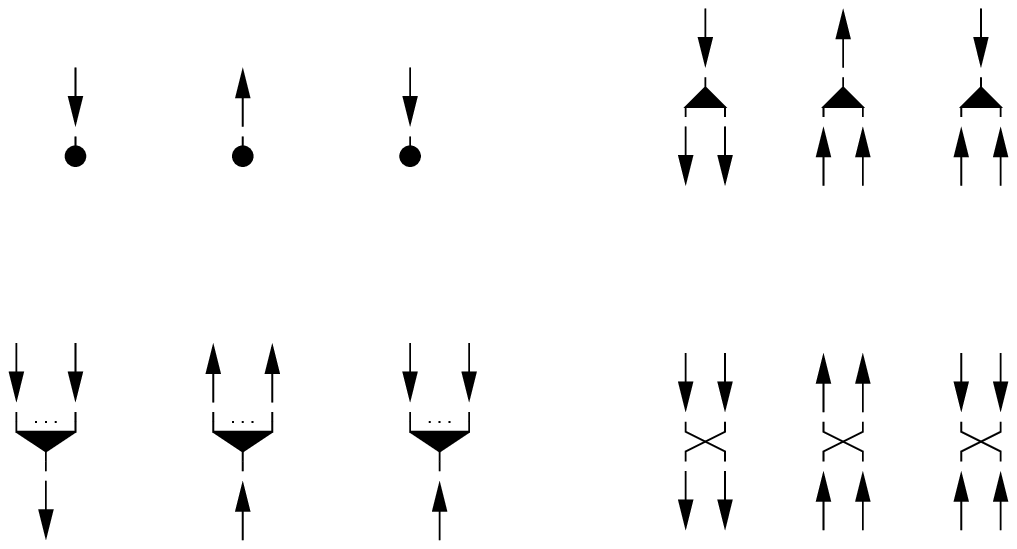}%
\end{picture}%
\setlength{\unitlength}{4144sp}%
\begingroup\makeatletter\ifx\SetFigFont\undefined%
\gdef\SetFigFont#1#2#3#4#5{%
  \reset@font\fontsize{#1}{#2pt}%
  \fontfamily{#3}\fontseries{#4}\fontshape{#5}%
  \selectfont}%
\fi\endgroup%
\begin{picture}(4888,2869)(153,-2054)
\put(2386,-106){\makebox(0,0)[lb]{\smash{{\SetFigFont{12}{14.4}{\rmdefault}{\mddefault}{\updefault}{\color[rgb]{0,0,0}$0$}%
}}}}
\put(3511,659){\makebox(0,0)[b]{\smash{{\SetFigFont{12}{14.4}{\rmdefault}{\mddefault}{\updefault}{\color[rgb]{0,0,0}$i$}%
}}}}
\put(4771,659){\makebox(0,0)[b]{\smash{{\SetFigFont{12}{14.4}{\rmdefault}{\mddefault}{\updefault}{\color[rgb]{0,0,0}$i$}%
}}}}
\put(631,389){\makebox(0,0)[b]{\smash{{\SetFigFont{12}{14.4}{\rmdefault}{\mddefault}{\updefault}{\color[rgb]{0,0,0}$i$}%
}}}}
\put(1396,389){\makebox(0,0)[b]{\smash{{\SetFigFont{12}{14.4}{\rmdefault}{\mddefault}{\updefault}{\color[rgb]{0,0,0}$1$}%
}}}}
\put(2161,389){\makebox(0,0)[b]{\smash{{\SetFigFont{12}{14.4}{\rmdefault}{\mddefault}{\updefault}{\color[rgb]{0,0,0}$i$}%
}}}}
\put(4141,659){\makebox(0,0)[b]{\smash{{\SetFigFont{12}{14.4}{\rmdefault}{\mddefault}{\updefault}{\color[rgb]{0,0,0}$i+j+1$}%
}}}}
\put(3421,-376){\makebox(0,0)[b]{\smash{{\SetFigFont{12}{14.4}{\rmdefault}{\mddefault}{\updefault}{\color[rgb]{0,0,0}$i$}%
}}}}
\put(3601,-376){\makebox(0,0)[b]{\smash{{\SetFigFont{12}{14.4}{\rmdefault}{\mddefault}{\updefault}{\color[rgb]{0,0,0}$i$}%
}}}}
\put(4051,-376){\makebox(0,0)[b]{\smash{{\SetFigFont{12}{14.4}{\rmdefault}{\mddefault}{\updefault}{\color[rgb]{0,0,0}$i$}%
}}}}
\put(4231,-376){\makebox(0,0)[b]{\smash{{\SetFigFont{12}{14.4}{\rmdefault}{\mddefault}{\updefault}{\color[rgb]{0,0,0}$j$}%
}}}}
\put(4681,-376){\makebox(0,0)[b]{\smash{{\SetFigFont{12}{14.4}{\rmdefault}{\mddefault}{\updefault}{\color[rgb]{0,0,0}$j$}%
}}}}
\put(4861,-376){\makebox(0,0)[b]{\smash{{\SetFigFont{12}{14.4}{\rmdefault}{\mddefault}{\updefault}{\color[rgb]{0,0,0}$k$}%
}}}}
\put(1261,-871){\makebox(0,0)[b]{\smash{{\SetFigFont{12}{14.4}{\rmdefault}{\mddefault}{\updefault}{\color[rgb]{0,0,0}$i$}%
}}}}
\put(1531,-871){\makebox(0,0)[b]{\smash{{\SetFigFont{12}{14.4}{\rmdefault}{\mddefault}{\updefault}{\color[rgb]{0,0,0}$i$}%
}}}}
\put(2161,-871){\makebox(0,0)[b]{\smash{{\SetFigFont{12}{14.4}{\rmdefault}{\mddefault}{\updefault}{\color[rgb]{0,0,0}$i_1$}%
}}}}
\put(2431,-871){\makebox(0,0)[b]{\smash{{\SetFigFont{12}{14.4}{\rmdefault}{\mddefault}{\updefault}{\color[rgb]{0,0,0}$i_n$}%
}}}}
\put(361,-871){\makebox(0,0)[b]{\smash{{\SetFigFont{12}{14.4}{\rmdefault}{\mddefault}{\updefault}{\color[rgb]{0,0,0}$i_1$}%
}}}}
\put(631,-871){\makebox(0,0)[b]{\smash{{\SetFigFont{12}{14.4}{\rmdefault}{\mddefault}{\updefault}{\color[rgb]{0,0,0}$i_n$}%
}}}}
\put(1396,-1996){\makebox(0,0)[b]{\smash{{\SetFigFont{12}{14.4}{\rmdefault}{\mddefault}{\updefault}{\color[rgb]{0,0,0}$i$}%
}}}}
\put(2296,-1996){\makebox(0,0)[b]{\smash{{\SetFigFont{12}{14.4}{\rmdefault}{\mddefault}{\updefault}{\color[rgb]{0,0,0}$j$}%
}}}}
\put(496,-1996){\makebox(0,0)[b]{\smash{{\SetFigFont{12}{14.4}{\rmdefault}{\mddefault}{\updefault}{\color[rgb]{0,0,0}$i_1+\dots+i_n+1$}%
}}}}
\put(3421,-916){\makebox(0,0)[b]{\smash{{\SetFigFont{12}{14.4}{\rmdefault}{\mddefault}{\updefault}{\color[rgb]{0,0,0}$i$}%
}}}}
\put(4051,-1951){\makebox(0,0)[b]{\smash{{\SetFigFont{12}{14.4}{\rmdefault}{\mddefault}{\updefault}{\color[rgb]{0,0,0}$i$}%
}}}}
\put(4231,-916){\makebox(0,0)[b]{\smash{{\SetFigFont{12}{14.4}{\rmdefault}{\mddefault}{\updefault}{\color[rgb]{0,0,0}$i$}%
}}}}
\put(4681,-916){\makebox(0,0)[b]{\smash{{\SetFigFont{12}{14.4}{\rmdefault}{\mddefault}{\updefault}{\color[rgb]{0,0,0}$i$}%
}}}}
\put(3421,-1951){\makebox(0,0)[b]{\smash{{\SetFigFont{12}{14.4}{\rmdefault}{\mddefault}{\updefault}{\color[rgb]{0,0,0}$j$}%
}}}}
\put(3601,-916){\makebox(0,0)[b]{\smash{{\SetFigFont{12}{14.4}{\rmdefault}{\mddefault}{\updefault}{\color[rgb]{0,0,0}$j$}%
}}}}
\put(4051,-916){\makebox(0,0)[b]{\smash{{\SetFigFont{12}{14.4}{\rmdefault}{\mddefault}{\updefault}{\color[rgb]{0,0,0}$j$}%
}}}}
\put(4231,-1951){\makebox(0,0)[b]{\smash{{\SetFigFont{12}{14.4}{\rmdefault}{\mddefault}{\updefault}{\color[rgb]{0,0,0}$j$}%
}}}}
\put(4861,-916){\makebox(0,0)[b]{\smash{{\SetFigFont{12}{14.4}{\rmdefault}{\mddefault}{\updefault}{\color[rgb]{0,0,0}$j$}%
}}}}
\put(4681,-1951){\makebox(0,0)[b]{\smash{{\SetFigFont{12}{14.4}{\rmdefault}{\mddefault}{\updefault}{\color[rgb]{0,0,0}$k$}%
}}}}
\put(4861,-1951){\makebox(0,0)[b]{\smash{{\SetFigFont{12}{14.4}{\rmdefault}{\mddefault}{\updefault}{\color[rgb]{0,0,0}$l$}%
}}}}
\put(3601,-1951){\makebox(0,0)[b]{\smash{{\SetFigFont{12}{14.4}{\rmdefault}{\mddefault}{\updefault}{\color[rgb]{0,0,0}$i$}%
}}}}
\put(5041,119){\makebox(0,0)[lb]{\smash{{\SetFigFont{12}{14.4}{\rmdefault}{\mddefault}{\updefault}{\color[rgb]{0,0,0}$\ul{i}+\ul{k}$}%
}}}}
\put(2611,-1411){\makebox(0,0)[lb]{\smash{{\SetFigFont{12}{14.4}{\rmdefault}{\mddefault}{\updefault}{\color[rgb]{0,0,0}$\ul{j}$}%
}}}}
\put(5041,-1411){\makebox(0,0)[lb]{\smash{{\SetFigFont{12}{14.4}{\rmdefault}{\mddefault}{\updefault}{\color[rgb]{0,0,0}$ij.\ul{l}+l.\ul{i+j}$}%
}}}}
\end{picture}%
\end{center}

\noindent Three diagrams are given for each operator $\phi$: two represent the functions $\phi_*$ and $\phi^*$ (how~$\phi$ transmits the current intensities) and one represents $[\phi]$ (the heat $\phi$ produces). Now, it is checked that, for every rule $f\fl g$ in $R_{\Delta\Sigma}$, the inequality $F(f)\succ F(g)$ holds, except for the rule $\alpha:s\alpha\fl t\alpha$, for which $F(s\alpha)=F(t\alpha)$. Let us check the (in)equalities for three sample rules. The complete computations are in [Guiraud 2004]. Let us start with the coassociativity rule for $\delta$:
\begin{center}
\begin{picture}(0,0)%
\includegraphics{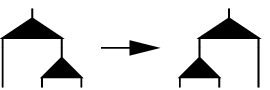}%
\end{picture}%
\setlength{\unitlength}{4144sp}%
\begingroup\makeatletter\ifx\SetFigFont\undefined%
\gdef\SetFigFont#1#2#3#4#5{%
  \reset@font\fontsize{#1}{#2pt}%
  \fontfamily{#3}\fontseries{#4}\fontshape{#5}%
  \selectfont}%
\fi\endgroup%
\begin{picture}(1194,384)(79,287)
\end{picture}%
\end{center}

\noindent One checks that the first two non-strict inequalities are satisfied:
$$\begin{cases}
((1\tens\delta)\circ\delta)_*(i)=(i,i,i)=((\delta\tens 1)\circ\delta)_*(i) \\
((1\tens\delta)\circ\delta)^*(i,j,k)=i+j+k+2=((\delta\tens 1)\circ\delta)^*(i,j,k).
\end{cases}$$

\noindent Moreover:
$\quad\begin{cases}
[(1\tens\delta)\circ\delta](i,j,k,l)=2.\ul{i}+\ul{l}+\ul{k+l+2} \\
[(\delta\tens 1)\circ\delta](i,j,k,l)=2.\ul{i}+\ul{l}+\ul{k}.
\end{cases}$ \\

\noindent Since $l+2>0$, one gets $\ul{k+l+2}>\ul{k}$ and the required strict inequality. Then, consider the rule $\alpha$ for which the chosen values do not work. One gets the two following equalities:
$$\begin{cases}
((1\tens\tau)\circ(\tau\tens 1)\circ(1\tens\tau))_*(i,j,k)=(k,j,i)=((\tau\tens 1)\circ(1\tens\tau)\circ(\tau\tens 1))_*(i,j,k) \\
((1\tens\tau)\circ(\tau\tens 1)\circ(1\tens\tau))^*(i,j,k)=(k,j,i)=((\tau\tens 1)\circ(1\tens\tau)\circ(\tau\tens 1))^*(i,j,k). \\
\end{cases}$$

\noindent And also this equality:
\begin{align*}
  & [(1\tens\tau)\circ(\tau\tens 1)\circ(1\tens\tau)](i,j,k,l,m,n) \\
=\quad & jk.\ul{m}+m.\ul{j+k}+i.(j+k).\ul{n}+n.(\ul{i+j}+\ul{j+k}) \\
=\quad & [(\tau\tens 1)\circ(1\tens\tau)\circ(\tau\tens 1)](i,j,k,l,m,n).
\end{align*}

\noindent To finish with our examples, let us consider the most complicated rule of this presentation, namely the local duplication rule:
\begin{center}
\begin{picture}(0,0)%
\includegraphics{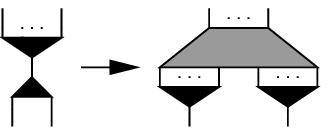}%
\end{picture}%
\setlength{\unitlength}{4144sp}%
\begingroup\makeatletter\ifx\SetFigFont\undefined%
\gdef\SetFigFont#1#2#3#4#5{%
  \reset@font\fontsize{#1}{#2pt}%
  \fontfamily{#3}\fontseries{#4}\fontshape{#5}%
  \selectfont}%
\fi\endgroup%
\begin{picture}(1464,564)(259,-793)
\end{picture}%
\end{center}

\noindent This is this rule that motivates the use of the rather complicated product category $\Or(\Nb^*,\Nb^*,[\Nb^*])$ to interpret $\mon{\Sigma^c}$. In order to make the computations for this rule, one must start by proving the following equations, which is done by iteration on the integer $n$:
$$\begin{cases}
(\delta_n)_*(i_1,\dots,i_n)=(i_1,\dots,i_n,i_1,\dots,i_n) \\
\delta_n^*(i_1,\dots,i_n,j_1,\dots,j_n)=(i_1+j_1+1,\dots,i_n+j_n+1) \\
[\delta_n](i_1,\dots,i_n,j_1,\dots,j_n,k_1,\dots,k_n) \\
\quad = \sum_{1\leq u\leq n}(\ul{i_u}+\ul{k_u})+\sum_{1\leq u<v\leq n}(i_ui_v.\ul{k_u}+k_u.\ul{i_u+i_v}).
\end{cases}$$

\noindent Then one gets these two equalities:
$$\begin{cases}
(\delta\circ\phi)_*(i_1,\dots,i_n)=(i_1+\dots+i_n+1,i_1+\dots+i_n+1)=((\phi\tens\phi)\circ\delta_n)_* \\
(\delta\circ\phi)^*(i,j)=(i+j+1,\dots,i+j+1)=((\phi\tens\phi)\circ\delta_n)^*.
\end{cases}$$

\noindent For the strict inequality to be checked:
$$\begin{cases}
[\delta\circ\phi](i_1,\dots,i_n,j,k)=\ul{j+k+1}+\ul{i_1+\dots+i_n+1}+\ul{k} \\
[(\phi\tens\phi)\circ\delta_n](i_1,\dots,i_n,j,k) \\
\quad = \ul{j}+(n+1+\sum_{1\leq u<v\leq n}i_ui_v).\ul{k}+\sum_{1\leq u<v\leq n}\ul{i_u}+k.\sum_{1\leq u<v\leq n}\ul{i_u+i_v}.
\end{cases}$$

\noindent The multiset order properties allow the conclusion: the left member of this rule is strictly greater than its right member. Indeed, it is a consequence from the following strict inequalities that hold in $[\Nb^*]$:
$$\begin{cases}
\ul{j+k+1}>\ul{j} \\
\ul{j+k+1}>\ul{k} \\
\ul{i_1+\dots+i_n+1}>\ul{i_u} &\text{for every }u \\
\ul{i_1+\dots+i_n+1}>\ul{i_u+i_v} &\text{for every }u\text{ and }v.
\end{cases}$$

\noindent The computations for the other rules are handled similarly, albeit more easily. Now, let us check the equivalence between termination of the $3$-polygraphs $(\Sigma^c,R_{\Delta\Sigma})$ and $(\Sigma^c,\{\alpha\})$. Since $\alpha$ is a rule of $R_{\Delta\Sigma}$, one concludes immediately that the termination of $(\Sigma^c,R_{\Delta\Sigma})$ implies the termination of $(\Sigma^c,\{\alpha\})$: any infinite reduction path generated by the latter would also be an infinite reduction path in the former. \\

\noindent Conversely, let us assume that $(\Sigma^c,\{\alpha\})$ terminates and that there exists an infinite reduction path $(f_n)_{n\in\Nb}$ in $(\Sigma^c,R_{\Delta\Sigma})$. This path yields an infinite decreasing sequence $(F(f_n))_n$ in $\Or(\Nb^*,\Nb^*,[\Nb^*])$, equipped with the order $\succeq$. Since this order terminates, the sequence is stationary, which means that there exists some natural number $n_0$ such that $F(f_n)=F(f_{n+1})$ whenever $n\geq n_0$. However, as proved earlier, one can have both $f\red{R_{\Delta\Sigma}}g$ and $F(f)=F(g)$ only if $f\red{\alpha}g$. This implies that the sequence $(f_n)_{n\geq n_0}$ is an infinite reduction path in $(\Sigma^c,\{\alpha\})$. However, the existence of such an infinite reduction path is prevented by the termination of $(\Sigma^c,\{\alpha\})$.
\findem\end{dem}

\noindent Now, there remains to prove that:

\begin{lem*}\label{lem2}
The $3$-polygraph $(\Sigma^c,\{\alpha\})$ terminates.
\end{lem*}

\begin{dem}
This is done using the technique from section \ref{technique}. The product category considered for the interpretations is $\Or(\Nb,\Nb,\Nb)$, where $\Nb$ is the set (or commutative monoid) of natural numbers, equipped with its natural order. We denote by $G$ the product category functor from $\mon{\Sigma^c}$ to $\Or(\Nb,\Nb,\Nb)$ defined by the following values on the operators of $\Sigma^c$:
\begin{center}
\begin{picture}(0,0)%
\includegraphics{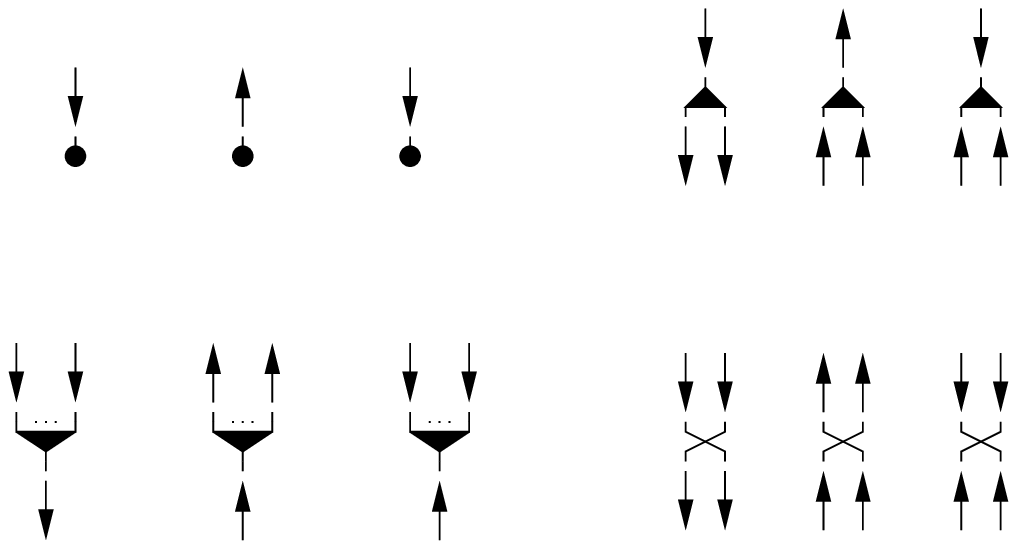}%
\end{picture}%
\setlength{\unitlength}{4144sp}%
\begingroup\makeatletter\ifx\SetFigFont\undefined%
\gdef\SetFigFont#1#2#3#4#5{%
  \reset@font\fontsize{#1}{#2pt}%
  \fontfamily{#3}\fontseries{#4}\fontshape{#5}%
  \selectfont}%
\fi\endgroup%
\begin{picture}(5272,2869)(153,-2054)
\put(3376,-1951){\makebox(0,0)[b]{\smash{{\SetFigFont{12}{14.4}{\rmdefault}{\mddefault}{\updefault}{\color[rgb]{0,0,0}$j$}%
}}}}
\put(631,389){\makebox(0,0)[b]{\smash{{\SetFigFont{12}{14.4}{\rmdefault}{\mddefault}{\updefault}{\color[rgb]{0,0,0}$i$}%
}}}}
\put(1396,389){\makebox(0,0)[b]{\smash{{\SetFigFont{12}{14.4}{\rmdefault}{\mddefault}{\updefault}{\color[rgb]{0,0,0}$1$}%
}}}}
\put(2161,389){\makebox(0,0)[b]{\smash{{\SetFigFont{12}{14.4}{\rmdefault}{\mddefault}{\updefault}{\color[rgb]{0,0,0}$i$}%
}}}}
\put(3511,659){\makebox(0,0)[b]{\smash{{\SetFigFont{12}{14.4}{\rmdefault}{\mddefault}{\updefault}{\color[rgb]{0,0,0}$i$}%
}}}}
\put(4771,659){\makebox(0,0)[b]{\smash{{\SetFigFont{12}{14.4}{\rmdefault}{\mddefault}{\updefault}{\color[rgb]{0,0,0}$i$}%
}}}}
\put(4141,659){\makebox(0,0)[b]{\smash{{\SetFigFont{12}{14.4}{\rmdefault}{\mddefault}{\updefault}{\color[rgb]{0,0,0}$i+j$}%
}}}}
\put(3421,-376){\makebox(0,0)[b]{\smash{{\SetFigFont{12}{14.4}{\rmdefault}{\mddefault}{\updefault}{\color[rgb]{0,0,0}$i$}%
}}}}
\put(3601,-376){\makebox(0,0)[b]{\smash{{\SetFigFont{12}{14.4}{\rmdefault}{\mddefault}{\updefault}{\color[rgb]{0,0,0}$i$}%
}}}}
\put(4051,-376){\makebox(0,0)[b]{\smash{{\SetFigFont{12}{14.4}{\rmdefault}{\mddefault}{\updefault}{\color[rgb]{0,0,0}$i$}%
}}}}
\put(4231,-376){\makebox(0,0)[b]{\smash{{\SetFigFont{12}{14.4}{\rmdefault}{\mddefault}{\updefault}{\color[rgb]{0,0,0}$j$}%
}}}}
\put(4681,-376){\makebox(0,0)[b]{\smash{{\SetFigFont{12}{14.4}{\rmdefault}{\mddefault}{\updefault}{\color[rgb]{0,0,0}$j$}%
}}}}
\put(4861,-376){\makebox(0,0)[b]{\smash{{\SetFigFont{12}{14.4}{\rmdefault}{\mddefault}{\updefault}{\color[rgb]{0,0,0}$k$}%
}}}}
\put(2161,-871){\makebox(0,0)[b]{\smash{{\SetFigFont{12}{14.4}{\rmdefault}{\mddefault}{\updefault}{\color[rgb]{0,0,0}$i_1$}%
}}}}
\put(1261,-871){\makebox(0,0)[b]{\smash{{\SetFigFont{12}{14.4}{\rmdefault}{\mddefault}{\updefault}{\color[rgb]{0,0,0}$i$}%
}}}}
\put(1531,-871){\makebox(0,0)[b]{\smash{{\SetFigFont{12}{14.4}{\rmdefault}{\mddefault}{\updefault}{\color[rgb]{0,0,0}$i$}%
}}}}
\put(1396,-1996){\makebox(0,0)[b]{\smash{{\SetFigFont{12}{14.4}{\rmdefault}{\mddefault}{\updefault}{\color[rgb]{0,0,0}$i$}%
}}}}
\put(2296,-1996){\makebox(0,0)[b]{\smash{{\SetFigFont{12}{14.4}{\rmdefault}{\mddefault}{\updefault}{\color[rgb]{0,0,0}$j$}%
}}}}
\put(3421,-916){\makebox(0,0)[b]{\smash{{\SetFigFont{12}{14.4}{\rmdefault}{\mddefault}{\updefault}{\color[rgb]{0,0,0}$i$}%
}}}}
\put(4231,-916){\makebox(0,0)[b]{\smash{{\SetFigFont{12}{14.4}{\rmdefault}{\mddefault}{\updefault}{\color[rgb]{0,0,0}$i$}%
}}}}
\put(4681,-916){\makebox(0,0)[b]{\smash{{\SetFigFont{12}{14.4}{\rmdefault}{\mddefault}{\updefault}{\color[rgb]{0,0,0}$i$}%
}}}}
\put(3601,-916){\makebox(0,0)[b]{\smash{{\SetFigFont{12}{14.4}{\rmdefault}{\mddefault}{\updefault}{\color[rgb]{0,0,0}$j$}%
}}}}
\put(4051,-916){\makebox(0,0)[b]{\smash{{\SetFigFont{12}{14.4}{\rmdefault}{\mddefault}{\updefault}{\color[rgb]{0,0,0}$j$}%
}}}}
\put(4861,-916){\makebox(0,0)[b]{\smash{{\SetFigFont{12}{14.4}{\rmdefault}{\mddefault}{\updefault}{\color[rgb]{0,0,0}$j$}%
}}}}
\put(4051,-1951){\makebox(0,0)[b]{\smash{{\SetFigFont{12}{14.4}{\rmdefault}{\mddefault}{\updefault}{\color[rgb]{0,0,0}$i$}%
}}}}
\put(4231,-1951){\makebox(0,0)[b]{\smash{{\SetFigFont{12}{14.4}{\rmdefault}{\mddefault}{\updefault}{\color[rgb]{0,0,0}$j$}%
}}}}
\put(4681,-1951){\makebox(0,0)[b]{\smash{{\SetFigFont{12}{14.4}{\rmdefault}{\mddefault}{\updefault}{\color[rgb]{0,0,0}$k$}%
}}}}
\put(4861,-1951){\makebox(0,0)[b]{\smash{{\SetFigFont{12}{14.4}{\rmdefault}{\mddefault}{\updefault}{\color[rgb]{0,0,0}$l$}%
}}}}
\put(631,-871){\makebox(0,0)[b]{\smash{{\SetFigFont{12}{14.4}{\rmdefault}{\mddefault}{\updefault}{\color[rgb]{0,0,0}$i_n$}%
}}}}
\put(361,-871){\makebox(0,0)[b]{\smash{{\SetFigFont{12}{14.4}{\rmdefault}{\mddefault}{\updefault}{\color[rgb]{0,0,0}$i_1$}%
}}}}
\put(2431,-871){\makebox(0,0)[b]{\smash{{\SetFigFont{12}{14.4}{\rmdefault}{\mddefault}{\updefault}{\color[rgb]{0,0,0}$i_n$}%
}}}}
\put(496,-1996){\makebox(0,0)[b]{\smash{{\SetFigFont{12}{14.4}{\rmdefault}{\mddefault}{\updefault}{\color[rgb]{0,0,0}$i_1+\dots+i_n$}%
}}}}
\put(2341,-106){\makebox(0,0)[lb]{\smash{{\SetFigFont{12}{14.4}{\rmdefault}{\mddefault}{\updefault}{\color[rgb]{0,0,0}$0$}%
}}}}
\put(2611,-1411){\makebox(0,0)[lb]{\smash{{\SetFigFont{12}{14.4}{\rmdefault}{\mddefault}{\updefault}{\color[rgb]{0,0,0}$0$}%
}}}}
\put(5041,-1411){\makebox(0,0)[lb]{\smash{{\SetFigFont{12}{14.4}{\rmdefault}{\mddefault}{\updefault}{\color[rgb]{0,0,0}$i+j$}%
}}}}
\put(5041,164){\makebox(0,0)[lb]{\smash{{\SetFigFont{12}{14.4}{\rmdefault}{\mddefault}{\updefault}{\color[rgb]{0,0,0}$0$}%
}}}}
\put(3691,-1951){\makebox(0,0)[b]{\smash{{\SetFigFont{12}{14.4}{\rmdefault}{\mddefault}{\updefault}{\color[rgb]{0,0,0}$i+1$}%
}}}}
\end{picture}%
\end{center}

\noindent We must check that $\alpha:s\alpha\fl t\alpha$ satisfies $F(s\alpha)\succ F(t\alpha)$. The computations give, on one hand, the two equalities:
$$\begin{cases}
((1\tens\tau)\circ(\tau\tens 1)\circ(1\tens\tau))_*(i,j,k)=(k,j+1,i+2)=((\tau\tens 1)\circ(1\tens\tau)\circ(\tau\tens 1))_*(i,j,k) \\
((1\tens\tau)\circ(\tau\tens 1)\circ(1\tens\tau))^*(i,j,k)=(k,j,i)=((\tau\tens 1)\circ(1\tens\tau)\circ(\tau\tens 1))^*(i,j,k). \\
\end{cases}$$

\noindent On the other hand, one gets:
$$\begin{cases}
[(1\tens\tau)\circ(\tau\tens 1)\circ(1\tens\tau)](i,j,k,l,m,n)=2i+2j+2k+2 \\
[(\tau\tens 1)\circ(1\tens\tau)\circ(\tau\tens 1)](i,j,k,l,m,n)=2i+2j+2k+1.
\end{cases}$$

\noindent By corollary \ref{cor2}, this gives the result.
\findem\end{dem}

\noindent Thus, one gets, as a corollary of lemmas \ref{lem1} and \ref{lem2}:

\begin{prop*}
The $3$-polygraph $(\Sigma^c,R_{\Delta\Sigma})$ terminates.
\end{prop*}

\noindent We recall the following result from [Guiraud 2004, proposition 5.31]:

\begin{prop*}
The $3$-polygraph $(\Sigma^c,R_{\Delta\Sigma})$ is locally confluent.
\end{prop*}

\noindent Finally, Newman's lemma [Baader Nipkow 1998] is applied to get theorem \ref{thm0}.

\subsection{Termination of $\mathbf{3}$-polygraph built from a terminating rewriting system}

\noindent This paragraph contains the proof of theorem \ref{prop1}, point 1: if a term rewriting system $(\Sigma,R)$ terminates, then so does its associated $3$-polygraph $(\Sigma^c,R^c)$. The proof once again uses a termination order obtained with theorem \ref{thm2}. However, integer values cannot be used here, since rules in $R$ are unknown. To handle this issue, the following classical result - see [Baader Nipkow 1998] - is used:

\begin{thm*}
A term rewriting system terminates if and only if there exists some mapping $|\cdot|$ from the set of terms $T\Sigma$ to $\Nb$ such that $|u|>|v|$ whenever $u$ is a term that reduces into another term $v$. Moreover, in that case, the mapping $|\cdot|$ can be chosen such that $|u|\geq|u'|$ whenever $u'$ is a subterm of $u$; the mapping can also be chosen so that it takes its values in any countable set.
\end{thm*}

\begin{dem}
If $(\Sigma,R)$ terminates, one can choose the mapping $|\cdot|$ to send each term $u$ onto the length of the longest reduction path starting from $u$; this mapping satisfies $|u|\geq|u'|$ if $u'$ is a subterm of $u$, since every reduction path from $u'$ yields a reduction path of the same length from $u$. Conversely, if such a mapping exists, an infinite reduction path $(u_n)_{n\in\Nb}$ in $(\Sigma,R)$ would generate a strictly decreasing infinite sequence $(|u_n|)_{n\in\Nb}$ in $\Nb$, which cannot exist; hence the term rewriting system $(\Sigma,R)$ terminates. If this is the case, the mapping $|\cdot|$ can be composed with any bijection $\sigma:\Nb\fl E$, where $E$ is any countable set.
\findem\end{dem}

\noindent Hence, from our terminating term rewriting system $(\Sigma,R)$, a mapping $|\cdot|:T\Sigma\fl\Nb^*$ is assumed to be chosen such that $|u|>|v|$ whenever $u$ reduces in $v$ and $|u|\geq|u'|$ whenever $u'$ is a subterm of $u$. From this mapping, one defines a binary relation $>$ on $T\Sigma$ by $u>v$ if, for every term context $c$, the inequality $|c[u]|>|c[v]|$ holds. From the fact that the usual order $>$ on $\Nb^*$ is a terminating strict order, this binary relation is proved to satisfy:

\begin{lem*}
The aforedefined binary relation $>$ on $T\Sigma$ is a terminating strict order.
\end{lem*}

\noindent Then, one builds the lexicographical order $\geq$ on $T\Sigma\times\Nb^*$: for this order, $(u,i)\geq(v,j)$ if $u>v$ or if $u=v$ and $i\geq j$. This order satisfies:

\begin{lem*}
This relation $\geq$ is an order on $T\Sigma\times\Nb^*$. Moreover, its strict part $>$ is a terminating strict order on $T\Sigma\times\Nb^*$.
\end{lem*}

\noindent The set $T\Sigma\times\Nb^*$, together with the aforedefined order, is taken as the first set used in the interpretation. The second one is a one-element set $\{\ast\}$ with the only possible order. Finally, the commutative monoid is once again $[\Nb^*]$ with its already-used multiset order. The product category~$\Or(T\Sigma\times\Nb^*,\{\ast\},[\Nb^*])$ is denoted by $\Or$.

\pagebreak
\noindent Sometimes, the two elements $((u_1,i_1),\dots,(u_n,i_n))$ of $(T\Sigma\times\Nb^*)^n$ and $(u_1,\dots,u_n;i_1,\dots,i_n)$ of~$T^n\times(\Nb^*)^n$ are identified.

The considered product category functor $F$ from $\mon{\Sigma^c}$ to $\Or$ is given by the following values (only two are given for each operator since the contravariant interpretation is trivial):
\begin{center}
\begin{picture}(0,0)%
\includegraphics{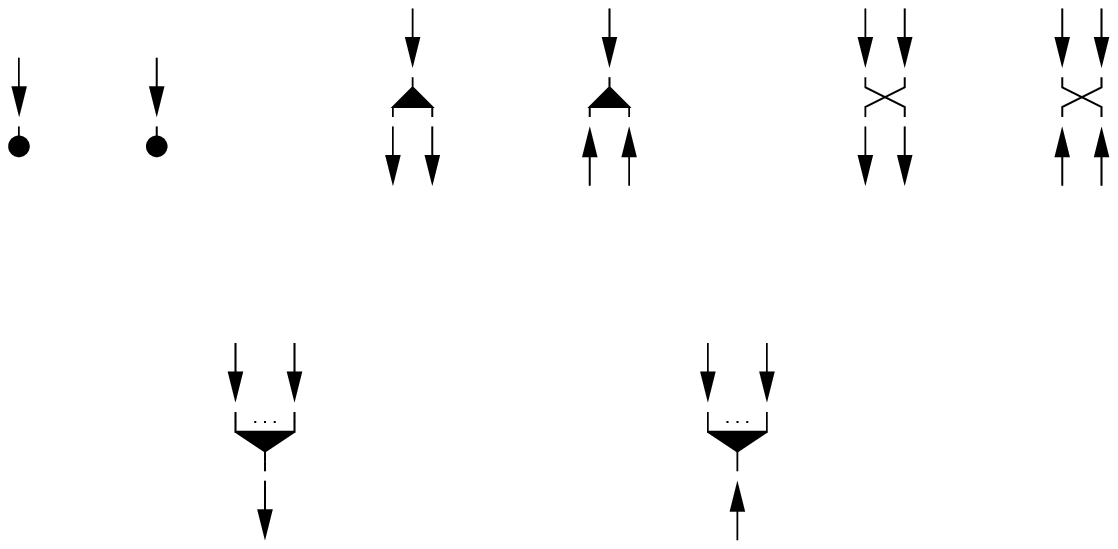}%
\end{picture}%
\setlength{\unitlength}{4144sp}%
\begingroup\makeatletter\ifx\SetFigFont\undefined%
\gdef\SetFigFont#1#2#3#4#5{%
  \reset@font\fontsize{#1}{#2pt}%
  \fontfamily{#3}\fontseries{#4}\fontshape{#5}%
  \selectfont}%
\fi\endgroup%
\begin{picture}(5667,2869)(-530,-2054)
\put(1261,-871){\makebox(0,0)[b]{\smash{{\SetFigFont{12}{14.4}{\rmdefault}{\mddefault}{\updefault}{\color[rgb]{0,0,0}$(u_n,i_n)$}%
}}}}
\put(-269,434){\makebox(0,0)[b]{\smash{{\SetFigFont{12}{14.4}{\rmdefault}{\mddefault}{\updefault}{\color[rgb]{0,0,0}$(u,i)$}%
}}}}
\put(361,434){\makebox(0,0)[b]{\smash{{\SetFigFont{12}{14.4}{\rmdefault}{\mddefault}{\updefault}{\color[rgb]{0,0,0}$(u,i)$}%
}}}}
\put(541,-61){\makebox(0,0)[lb]{\smash{{\SetFigFont{12}{14.4}{\rmdefault}{\mddefault}{\updefault}{\color[rgb]{0,0,0}$0$}%
}}}}
\put(1531,659){\makebox(0,0)[b]{\smash{{\SetFigFont{12}{14.4}{\rmdefault}{\mddefault}{\updefault}{\color[rgb]{0,0,0}$(u,i)$}%
}}}}
\put(1306,-376){\makebox(0,0)[b]{\smash{{\SetFigFont{12}{14.4}{\rmdefault}{\mddefault}{\updefault}{\color[rgb]{0,0,0}$(u,i)$}%
}}}}
\put(1756,-376){\makebox(0,0)[b]{\smash{{\SetFigFont{12}{14.4}{\rmdefault}{\mddefault}{\updefault}{\color[rgb]{0,0,0}$(u,i)$}%
}}}}
\put(2431,659){\makebox(0,0)[b]{\smash{{\SetFigFont{12}{14.4}{\rmdefault}{\mddefault}{\updefault}{\color[rgb]{0,0,0}$(u,i)$}%
}}}}
\put(2521,-376){\makebox(0,0)[b]{\smash{{\SetFigFont{12}{14.4}{\rmdefault}{\mddefault}{\updefault}{\color[rgb]{0,0,0}$\ast$}%
}}}}
\put(2701,164){\makebox(0,0)[lb]{\smash{{\SetFigFont{12}{14.4}{\rmdefault}{\mddefault}{\updefault}{\color[rgb]{0,0,0}$i.\ul{u}$}%
}}}}
\put(2341,-376){\makebox(0,0)[b]{\smash{{\SetFigFont{12}{14.4}{\rmdefault}{\mddefault}{\updefault}{\color[rgb]{0,0,0}$\ast$}%
}}}}
\put(3466,659){\makebox(0,0)[b]{\smash{{\SetFigFont{12}{14.4}{\rmdefault}{\mddefault}{\updefault}{\color[rgb]{0,0,0}$(u,i)$}%
}}}}
\put(3466,-376){\makebox(0,0)[b]{\smash{{\SetFigFont{12}{14.4}{\rmdefault}{\mddefault}{\updefault}{\color[rgb]{0,0,0}$(v,j)$}%
}}}}
\put(3916,-376){\makebox(0,0)[b]{\smash{{\SetFigFont{12}{14.4}{\rmdefault}{\mddefault}{\updefault}{\color[rgb]{0,0,0}$(u,i)$}%
}}}}
\put(3916,659){\makebox(0,0)[b]{\smash{{\SetFigFont{12}{14.4}{\rmdefault}{\mddefault}{\updefault}{\color[rgb]{0,0,0}$(v,j)$}%
}}}}
\put(4366,659){\makebox(0,0)[b]{\smash{{\SetFigFont{12}{14.4}{\rmdefault}{\mddefault}{\updefault}{\color[rgb]{0,0,0}$(u,i)$}%
}}}}
\put(4501,-376){\makebox(0,0)[b]{\smash{{\SetFigFont{12}{14.4}{\rmdefault}{\mddefault}{\updefault}{\color[rgb]{0,0,0}$\ast$}%
}}}}
\put(4681,-376){\makebox(0,0)[b]{\smash{{\SetFigFont{12}{14.4}{\rmdefault}{\mddefault}{\updefault}{\color[rgb]{0,0,0}$\ast$}%
}}}}
\put(4861,164){\makebox(0,0)[lb]{\smash{{\SetFigFont{12}{14.4}{\rmdefault}{\mddefault}{\updefault}{\color[rgb]{0,0,0}$0$}%
}}}}
\put(4816,659){\makebox(0,0)[b]{\smash{{\SetFigFont{12}{14.4}{\rmdefault}{\mddefault}{\updefault}{\color[rgb]{0,0,0}$(v,j)$}%
}}}}
\put(3016,-1996){\makebox(0,0)[b]{\smash{{\SetFigFont{12}{14.4}{\rmdefault}{\mddefault}{\updefault}{\color[rgb]{0,0,0}$\ast$}%
}}}}
\put(3331,-1411){\makebox(0,0)[lb]{\smash{{\SetFigFont{12}{14.4}{\rmdefault}{\mddefault}{\updefault}{\color[rgb]{0,0,0}$(i_1+\dots+i_n).\ul{|\phi(u_1,\dots,u_n)|}$}%
}}}}
\put(2656,-871){\makebox(0,0)[b]{\smash{{\SetFigFont{12}{14.4}{\rmdefault}{\mddefault}{\updefault}{\color[rgb]{0,0,0}$(u_1,i_1)$}%
}}}}
\put(3421,-871){\makebox(0,0)[b]{\smash{{\SetFigFont{12}{14.4}{\rmdefault}{\mddefault}{\updefault}{\color[rgb]{0,0,0}$(u_n,i_n)$}%
}}}}
\put(856,-1996){\makebox(0,0)[b]{\smash{{\SetFigFont{12}{14.4}{\rmdefault}{\mddefault}{\updefault}{\color[rgb]{0,0,0}$(\phi(u_1,\dots,u_n),2.(i_1+\dots+i_n))$}%
}}}}
\put(496,-871){\makebox(0,0)[b]{\smash{{\SetFigFont{12}{14.4}{\rmdefault}{\mddefault}{\updefault}{\color[rgb]{0,0,0}$(u_1,i_1)$}%
}}}}
\end{picture}%
\end{center}

\noindent There are two steps to check the conditions given in corollary \ref{cor2}: the first one consists in ensuring that each given operation is monotone; the second part is about computing if $F(f)>F(g)$ holds for every rule $f\fl g$ in $R^c$. \\

\noindent For the first part, consider, for example, the functions $\phi_*$ and $[\phi]$ for some fixed operator $\phi$ in $\Sigma(n,1)$, $n\geq 1$. Let us consider terms $u_1,\dots,u_n$, $v_1,\dots,v_n$ and non-zero integers $i_1,\dots,i_n$, $j_1,\dots,j_n$. Let us assume that $(u_k,i_k)\geq (v_k,j_k)$ for every $k$. In order to prove that $\phi_*$ is monotone, one must check that either $\phi(u_1,\dots,u_n)>\phi(v_1,\dots,v_n)$ or both are equal and $i_1+\dots+i_n\geq j_1+\dots+j_n$. Let $c$ be a context. Since, for every $k$, $u_k\geq v_k$ and $c\circ\phi(v_1,\dots,v_{k-1},\sq,u_{k+1},\dots,u_n)$ is a context, one gets the following inequality:
$$
|c[\phi(v_1,\dots,v_{k-1},u_k,\dots,u_n)]|\geq |c[\phi(v_1,\dots,v_k,u_{k+1},\dots,u_n)]|.
$$

\noindent Furthermore, if $u_k>v_k$ for some $k$, then this inequality is strict for the same $k$; in this case:
$$
|c[\phi(u_1,\dots,u_n)]|>|c[\phi(v_1,\dots,v_n)]|.
$$

\noindent Consequently, $\phi(u_1,\dots,u_n)>\phi(v_1,\dots,v_n)$. Otherwise, if $u_k=v_k$ and $i_k\geq j_k$ for every $k$, then:
$$\begin{cases}
|c[\phi(u_1,\dots,u_n)]|=|c[\phi(v_1,\dots,v_n)]| \\
i_1+\dots+i_n\geq j_1+\dots+j_n.
\end{cases}$$

\noindent Thus, in both cases:
$$
(\phi(u_1,\dots,u_n),2.(i_1+\dots+i_n))\geq (\phi(v_1,\dots,v_n),2.(j_1+\dots+j_n)).
$$

\noindent In order to prove that $[\phi]$ is monotone, let us fix some $k$ in $[n]$. Then, either $u_k>v_k$ or $u_k=v_k$ and $i_k\geq j_k$. In the first case:
$$
|\phi(v_1,\dots,v_{k-1},u_k,\dots,u_n)|>|\phi(v_1,\dots,v_k,u_{k+1},\dots,u_n)|.
$$

\noindent Thus, by definition of the multiset order on $[\Nb^*]$:
\begin{align*}
 & (j_1+\dots+j_{k-1}+i_k+\dots+i_n).\ul{|\phi(v_1,\dots,v_{k-1},u_k,\dots,u_n)|} \\
>\quad & (j_1+\dots+j_k+i_{k+1}+\dots+i_n).\ul{|\phi(v_1,\dots,v_k,u_{k+1},\dots,u_n)|}.
\end{align*}

\noindent In the second case, where $u_k=v_k$ and $i_k\geq j_k$:
\begin{align*}
& (j_1+\dots+j_{k-1}+i_k+\dots+i_n).\ul{|\phi(v_1,\dots,v_{k-1},u_k,\dots,u_n)|} \\
\geq \quad & (j_1+\dots+j_k+i_{k+1}+\dots+i_n).\ul{|\phi(v_1,\dots,v_k,u_{k+1},\dots,u_n)|}.
\end{align*}

\noindent Finally:
$$
(i_1+\dots+i_n).\ul{|\phi(u_1,\dots,u_n)|}\geq (j_1+\dots+j_n).\ul{|\phi(v_1,\dots,v_n)|}.
$$

\noindent If $\gamma$ is a constant in $\Sigma(0,1)$ or for operators in $\Delta$, proofs are direct. Furthermore, for each operator $\phi$ in either $\Sigma$ or $\Delta$, the operation $\phi^*$ is the only map from $\{\ast\}$ to itself, and it is monotone, so that:

\begin{lem*}
For every operator $\phi$ in $\Sigma^c$, the aforegiven functions $\phi_*$, $\phi^*$ and $[\phi]$ are monotone.
\end{lem*}

\noindent Then, we must check if $F(f)\succ F(g)$ for every rule $f\fl g$ in $R^c$. Let us recall that this family of rules consists of three subfamilies: $R_{\Delta}$, $R_{\Sigma}$ and $\Phi(R)$. For any rule $f\fl g$ in the first family $R_{\Delta}$, one gets $F(f)=F(g)$, except for left and right counit rules, where $F(f)\succ F(g)$. Computations for rules in $R_{\Sigma}$ are more complicated; let us examine, for example, the rule for local duplication and one of the rules for local permutation:
\begin{center}
\begin{picture}(0,0)%
\includegraphics{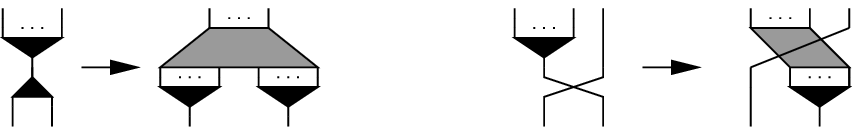}%
\end{picture}%
\setlength{\unitlength}{4144sp}%
\begingroup\makeatletter\ifx\SetFigFont\undefined%
\gdef\SetFigFont#1#2#3#4#5{%
  \reset@font\fontsize{#1}{#2pt}%
  \fontfamily{#3}\fontseries{#4}\fontshape{#5}%
  \selectfont}%
\fi\endgroup%
\begin{picture}(3894,564)(-11,107)
\end{picture}%
\end{center}

\noindent Let us fix some natural number $n\geq 1$ and some $\phi$ in $\Sigma(n,1)$; for constants in $\Sigma(0,1)$, computations are direct. By iteration on $n$, the following equalities are proved:
$$
(\delta_n)_*(\vec{u},\vec{\iota})=(\vec{u},\vec{u};\vec{\iota},\vec{\iota}) \quad\text{and}\quad [\delta_n](\vec{u};\vec{\iota})=i_1.\ul{u_1}+\dots+i_n.\ul{u_n}.
$$

\noindent This gives, at first:
\begin{align*}
(\delta\circ\phi)_*(\vec{u},\vec{\iota})
&= (\phi(\vec{u}),\phi(\vec{u});2.(i_1+\dots+i_n),2.(i_1+\dots+i_n)) \\
&= ((\phi\tens\phi)\circ\delta_n)_*(\vec{u},\vec{\iota}).
\end{align*}

\noindent Then: $[\delta\circ\phi](\vec{u},\vec{\iota})=3.(i_1+\dots+i_n).\ul{|\phi(\vec{u})|}$. To be compared with:
$$
[(\phi\tens\phi)\circ\delta_n](\vec{u},\vec{\iota})=2.(i_1+\dots+i_n).\ul{|\phi(\vec{u})|}+i_1.\ul{u_1}+\dots+i_n.\ul{u_n}.
$$

\noindent Since $u_k$ is a subterm of $\phi(\vec{u})$ for every $k$, and by asumption on $|\cdot|$, the inequality $|\phi(\vec{u})|\geq |u_k|$ holds. Hence, for every $k$:
$$
i_k.\ul{|\phi(\vec{u})|}\geq i_k.\ul{u_k}.
$$

\noindent Finally: $(i_1+\dots+i_n).\ul{|\phi(\vec{u})|}\geq i_1.\ul{u_1}+\dots+i_n.\ul{u_n}$. This gives the inequality $[\delta\circ\phi]\geq[(\phi\tens\phi)\circ\delta_n]$. Now, let us consider the first rule for local permutation; the first step is to prove, by iteration on $n$:
$$
(\tau_{n,1})_*(\vec{u},v;\vec{\iota},j)=(v,\vec{u};j,\vec{\iota}) \quad\text{and}\quad [\tau_{n,1}](\vec{u},v;\vec{\iota},j)=0.
$$

\noindent Then: $\quad (\tau\circ(\phi\tens 1))_*(\vec{u},v;\vec{i},j)\:=\:(v,\phi(\vec{u});j,2.(i_1+\dots+i_n))\:=\:((1\tens\phi)\circ\tau_{n,1})_*(\vec{u},v;\vec{i},j)$. \\

\noindent And: $\quad [\tau\circ(\phi\tens 1)](\vec{u},v;\vec{i},j)\:=\:(i_1+\dots+i_n).\ul{|\phi(\vec{u})|}\:=\:[(1\tens\phi)\circ\tau_{n,1}](\vec{u},v;\vec{i},j)$. \\

\noindent The other rules in $R_{\Sigma}$ are similarly handled and give similar results: for every rule $f\fl g$ in~$R_{\Sigma}$, the inequality $F(f)\succeq F(g)$ holds in $\Or$. The final part concerns the family $\Phi(R)$ of rules. Let us assume that $\alpha:f\fl g$ is a rule in $R$; its translation by $\Phi$ is the rule $\Phi(\alpha):\Phi(f)\fl\Phi^{\sharp f}(g)$. Let us prove that $F\circ\Phi(f)\succ F\circ\Phi^{\sharp f}(g)$. The first step is to prove, by iteration on the degree of terms in $T\Sigma$, the following lemma:

\begin{lem*}
Let $u$ be a term in $T\Sigma$, $n$ be an integer such that $n\geq \sharp u$, $\vec{v}$ a family of $n$ terms in $T\Sigma$ and~$\vec{\iota}$ a family of $n$ non-zero natural numbers. Let us denote by $\sigma_{\vec{v}}$ the substitution defined by $x_k\cdot\sigma_v=v_k$ if $k\leq n$ and $x_k$ otherwise. Then:
\begin{enumerate}
\item There exists some non-zero integer $k$ such that $(\Phi^n(u))_*(u\cdot\sigma_{\vec{v}},k)$.
\item The inequality $[\Phi^n(u)](\vec{v},\vec{\iota})<\ul{|u\cdot\sigma_{\vec{v}}+1|}$ holds in $[\Nb^*]$.
\item If $u$ is not a variable, then the inequality $[\Phi^n(u)](\vec{v},\vec{\iota})\geq \ul{|u\cdot\sigma_{\vec{v}}|}$ also holds in $[\Nb^*]$.
\end{enumerate}
\end{lem*}

\noindent Point 1 gives, when applied to $f$ and $g$ with $n=\sharp f$, the existence of non-zero natural numbers $k$ and~$k'$ such that $\Phi(f)_*(\vec{u},\vec{\iota})=(f\cdot\sigma_{\vec{u}},k)$ and $\Phi^n(g)_*(\vec{u},\vec{\iota})=(g\cdot\sigma_{\vec{u}},k')$. Let us consider some context $c$. By definition of the reduction relation $\red{\alpha}\:$ generated by the rule $\alpha$, one gets $c[f\cdot\sigma_{\vec{u}}]\red{\alpha}\:c[g\cdot\sigma_{\vec{u}}]$. Consequently, the properties of $|\cdot|$ give $|c[f\cdot\sigma_{\vec{u}}]|>|c[g\cdot\sigma_{\vec{u}}]|$. This holds for any context thus, by definition of $>$ on $T\Sigma$, one gets $f\cdot\sigma_{\vec{u}}>g\cdot\sigma_{\vec{u}}$. Finally, using the definition of $>$ on $T\Sigma\times\Nb^*$:
$$
\Phi_*(f)>\Phi_*(g).
$$

\noindent Let us prove now that $[\Phi(f)]>[\Phi^n(g)]$. Since $\alpha$ is a term rewrite rule, its source $f$ is a non-variable term. Hence, point 3 of the previous lemma gives the inequality $[\Phi(f)](\vec{u},\vec{\iota})\geq \ul{|f\cdot\sigma_{\vec{u}}|}$. Moreover, point~2 gives $[\Phi^n(g)](\vec{u},\vec{\iota})<\ul{|g\cdot\sigma_{\vec{u}}+1|}$. Finally, since the reduction $f\cdot\sigma_{\vec{u}}\red{\alpha}\:g\cdot\sigma_{\vec{u}}$ holds in $(\Sigma,R)$ and by properties of $|\cdot|$: $|f\cdot\sigma_{\vec{u}}|>|g\cdot\sigma_{\vec{v}}|$. There remains to concatenate these three inequalities to get $[\Phi(f)]>[\Phi^n(g)]$ and, as a consequence $F\circ\Phi(f)\succ F\circ\Phi^n(g)$. The product category functor $F$ from~$\mon{\Sigma^c}$ to $\Or$ gives us $F(f)\succ F(g)$ for every rule $f\fl g$ in $\Phi(R)$ and $F(f)\succeq F(g)$ for every rule $f\fl g$ in $R_{\Delta\Sigma}$. This yields the following result:

\begin{prop*}
If the term rewriting system $(\Sigma,R)$ terminates, then termination of the $3$-polygraph $(\Sigma^c,R^c)$ is equivalent to termination of $(\Sigma^c,R_{\Delta\Sigma})$.
\end{prop*}

\noindent Since we already know that $(\Sigma^c,R_{\Delta\Sigma})$ always terminates, this concludes the proof of theorem~\ref{thm1}.

\section{Application 2: a convergent $\mathbf{3}$-polygraph for a commutative equational theory}\label{app2}

\noindent This final section is devoted to give a convergent presentation of the equational theory of $\zb$-vector spaces, which is, as mentionned before, a commutative equational theory and thus do not have any convergent presentation by a term rewriting system. \\

\noindent In section \ref{equational}, we have considered three term rewriting systems $(\Sigma,R_0)$, $(\Sigma,R_1)$ and $(\Sigma,R_2)$ that repectively present the equational theories of monoids, of commutative monoids and of $\zb$-vector spaces. All three have two operators, a product and a unit, and they have respectively three, four and five rules. Thus, their associated $3$-polygraphs have five operators together with twenty-three rules for $(\Sigma^c,R_0^c)$, twenty-four for $(\Sigma^c,R_1^c)$ and twenty-five for $(\Sigma^c,R_2^c)$.

Since $(\Sigma,R_0)$ is a left-linear convergent term rewriting system, theorem \ref{thm1} ensures, in particular, that $(\Sigma^c,R_0^c)$ is a convergent presentation of the theory of monoids, \emph{with explicit resource management}. The term rewriting system $(\Sigma,R_1)$ is left-linear, non-terminating (due to the commutativity rule) and non-confluent (though it could be completed to get a confluent rewriting system), hence theorem \ref{thm1} gives us that $(\Sigma^c,R_1^c)$ is a non-terminating and non-confluent presentation of the equational theory of commutative monoids, with explicit resource management. Finally, the term rewriting system $(\Sigma,R_2)$ is a non-left-linear, non-terminating and non-confluent term rewriting system: non-left-linearity denies us any information coming from theorem \ref{thm1} about this presentation. \\

\noindent However, there is, in [Lafont 2003], an equivalent $3$-polygraph called $\lz$. Its signature contains a sixth operator, called $\kappa$ and pictured this way:
\begin{center}
\begin{picture}(0,0)%
\includegraphics{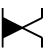}%
\end{picture}%
\setlength{\unitlength}{4144sp}%
\begingroup\makeatletter\ifx\SetFigFont\undefined%
\gdef\SetFigFont#1#2#3#4#5{%
  \reset@font\fontsize{#1}{#2pt}%
  \fontfamily{#3}\fontseries{#4}\fontshape{#5}%
  \selectfont}%
\fi\endgroup%
\begin{picture}(204,204)(439,467)
\end{picture}%
\end{center}

\noindent This new operator is said to be \emph{superfluous} since it represents, in a $\zb$-vector space, the concrete operation $\kappa(x,y)=(\mu(x,y),x)$ that can be expressed in terms of $\mu$, $\delta$ and $\tau$. In the presentation, this relation is enforced by means of the following extra rule:
\begin{center}
\begin{picture}(0,0)%
\includegraphics{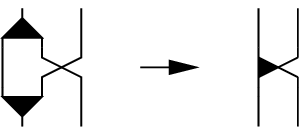}%
\end{picture}%
\setlength{\unitlength}{4144sp}%
\begingroup\makeatletter\ifx\SetFigFont\undefined%
\gdef\SetFigFont#1#2#3#4#5{%
  \reset@font\fontsize{#1}{#2pt}%
  \fontfamily{#3}\fontseries{#4}\fontshape{#5}%
  \selectfont}%
\fi\endgroup%
\begin{picture}(1374,564)(169,107)
\end{picture}%
\end{center}

\noindent The main objective of these new operator and rule is to make proof of termination easier (if not just possible). Then, one has to add a certain amount of rules in order to complete the presentation, to finally obtain the $3$-polygraph $\lz$, discovered and baptized in [Lafont 2003]. \\

\noindent This polygraph has six operators:
\begin{center}
\begin{picture}(0,0)%
\includegraphics{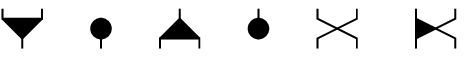}%
\end{picture}%
\setlength{\unitlength}{4144sp}%
\begingroup\makeatletter\ifx\SetFigFont\undefined%
\gdef\SetFigFont#1#2#3#4#5{%
  \reset@font\fontsize{#1}{#2pt}%
  \fontfamily{#3}\fontseries{#4}\fontshape{#5}%
  \selectfont}%
\fi\endgroup%
\begin{picture}(2094,520)(169,61)
\put(2161,119){\makebox(0,0)[b]{\smash{{\SetFigFont{12}{14.4}{\rmdefault}{\mddefault}{\updefault}$\kappa$}}}}
\put(271,119){\makebox(0,0)[b]{\smash{{\SetFigFont{12}{14.4}{\rmdefault}{\mddefault}{\updefault}$\mu$}}}}
\put(631,119){\makebox(0,0)[b]{\smash{{\SetFigFont{12}{14.4}{\rmdefault}{\mddefault}{\updefault}$\eta$}}}}
\put(991,119){\makebox(0,0)[b]{\smash{{\SetFigFont{12}{14.4}{\rmdefault}{\mddefault}{\updefault}$\delta$}}}}
\put(1351,119){\makebox(0,0)[b]{\smash{{\SetFigFont{12}{14.4}{\rmdefault}{\mddefault}{\updefault}$\epsilon$}}}}
\put(1711,119){\makebox(0,0)[b]{\smash{{\SetFigFont{12}{14.4}{\rmdefault}{\mddefault}{\updefault}$\tau$}}}}
\end{picture}%
\end{center}

\pagebreak
\noindent And sixty-seven rules:
\begin{center}
\includegraphics[width=16cm]{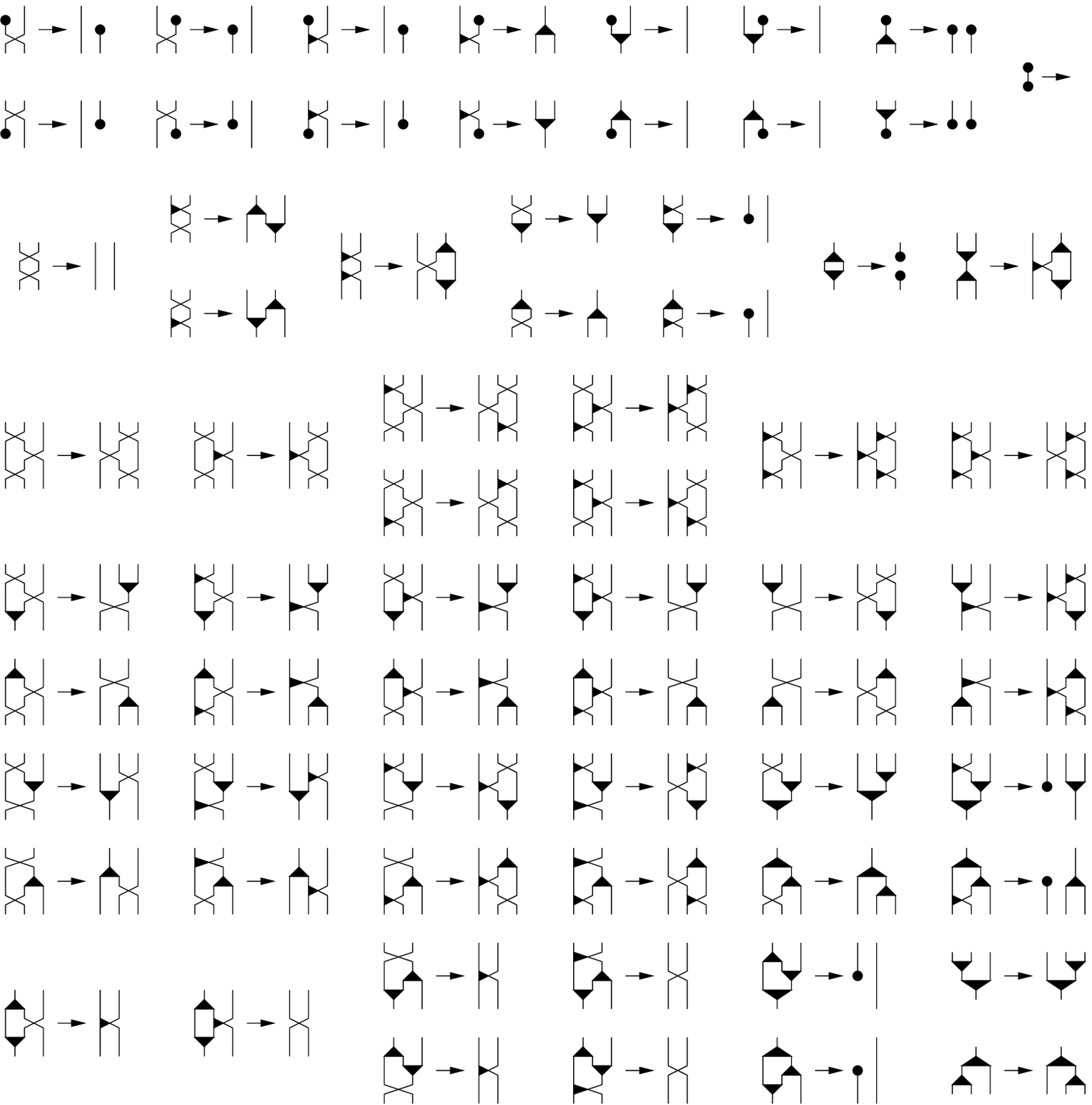}
\end{center}

\noindent From [Lafont 2003], we already know that this presentation is confluent but termination was still a conjecture. The technique presented in section \ref{technique} now allows us to prove that it is also terminating, hence convergent. The interpretation product category we use is $\Or(\Nb^*,\Nb^*,[\Nb^*])$, once again denoted by $\Or$. The interpretation functor $F$ is given by the following values on generating operators:
\begin{center}
\begin{picture}(0,0)%
\includegraphics{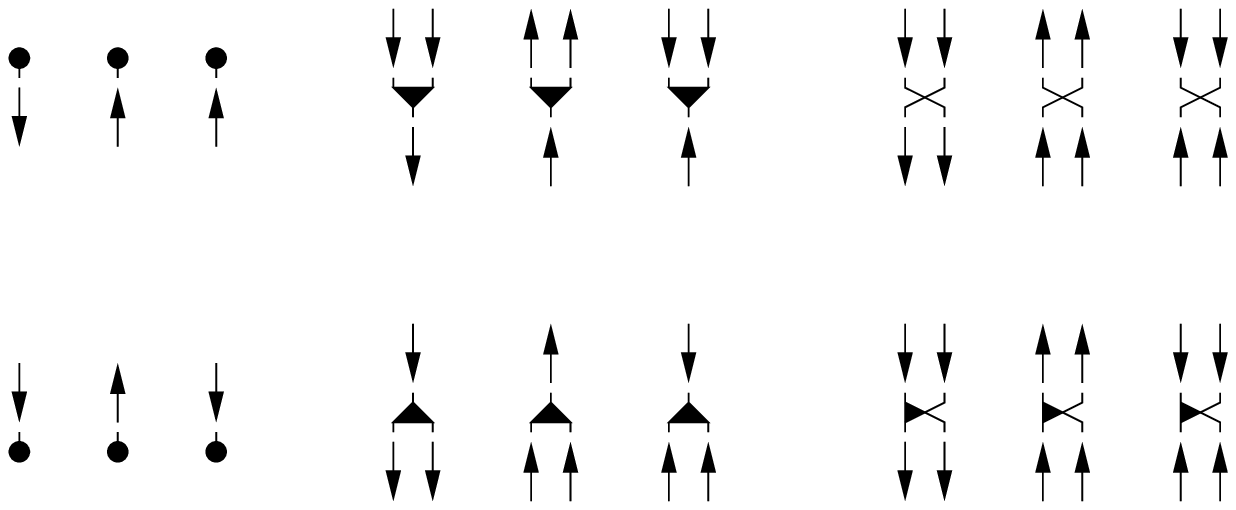}%
\end{picture}%
\setlength{\unitlength}{4144sp}%
\begingroup\makeatletter\ifx\SetFigFont\undefined%
\gdef\SetFigFont#1#2#3#4#5{%
  \reset@font\fontsize{#1}{#2pt}%
  \fontfamily{#3}\fontseries{#4}\fontshape{#5}%
  \selectfont}%
\fi\endgroup%
\begin{picture}(5809,2689)(42,-1784)
\put(5086,-691){\makebox(0,0)[b]{\smash{{\SetFigFont{12}{14.4}{\rmdefault}{\mddefault}{\updefault}$i$}}}}
\put(181,-106){\makebox(0,0)[b]{\smash{{\SetFigFont{12}{14.4}{\rmdefault}{\mddefault}{\updefault}$1$}}}}
\put(631,-871){\makebox(0,0)[b]{\smash{{\SetFigFont{12}{14.4}{\rmdefault}{\mddefault}{\updefault}$1$}}}}
\put(181,-871){\makebox(0,0)[b]{\smash{{\SetFigFont{12}{14.4}{\rmdefault}{\mddefault}{\updefault}$i$}}}}
\put(631,-106){\makebox(0,0)[b]{\smash{{\SetFigFont{12}{14.4}{\rmdefault}{\mddefault}{\updefault}$i$}}}}
\put(1081,-106){\makebox(0,0)[b]{\smash{{\SetFigFont{12}{14.4}{\rmdefault}{\mddefault}{\updefault}$i$}}}}
\put(1081,-871){\makebox(0,0)[b]{\smash{{\SetFigFont{12}{14.4}{\rmdefault}{\mddefault}{\updefault}$i$}}}}
\put(1261,389){\makebox(0,0)[lb]{\smash{{\SetFigFont{12}{14.4}{\rmdefault}{\mddefault}{\updefault}$\ul{i}$}}}}
\put(1261,-1366){\makebox(0,0)[lb]{\smash{{\SetFigFont{12}{14.4}{\rmdefault}{\mddefault}{\updefault}$\ul{i}$}}}}
\put(1891,749){\makebox(0,0)[b]{\smash{{\SetFigFont{12}{14.4}{\rmdefault}{\mddefault}{\updefault}$i$}}}}
\put(2521,749){\makebox(0,0)[b]{\smash{{\SetFigFont{12}{14.4}{\rmdefault}{\mddefault}{\updefault}$i$}}}}
\put(2701,749){\makebox(0,0)[b]{\smash{{\SetFigFont{12}{14.4}{\rmdefault}{\mddefault}{\updefault}$i$}}}}
\put(2611,-286){\makebox(0,0)[b]{\smash{{\SetFigFont{12}{14.4}{\rmdefault}{\mddefault}{\updefault}$i$}}}}
\put(3151,749){\makebox(0,0)[b]{\smash{{\SetFigFont{12}{14.4}{\rmdefault}{\mddefault}{\updefault}$i$}}}}
\put(1981,-691){\makebox(0,0)[b]{\smash{{\SetFigFont{12}{14.4}{\rmdefault}{\mddefault}{\updefault}$i$}}}}
\put(2521,-1726){\makebox(0,0)[b]{\smash{{\SetFigFont{12}{14.4}{\rmdefault}{\mddefault}{\updefault}$i$}}}}
\put(1891,-1726){\makebox(0,0)[b]{\smash{{\SetFigFont{12}{14.4}{\rmdefault}{\mddefault}{\updefault}$i$}}}}
\put(2071,-1726){\makebox(0,0)[b]{\smash{{\SetFigFont{12}{14.4}{\rmdefault}{\mddefault}{\updefault}$i$}}}}
\put(3241,-691){\makebox(0,0)[b]{\smash{{\SetFigFont{12}{14.4}{\rmdefault}{\mddefault}{\updefault}$i$}}}}
\put(2071,749){\makebox(0,0)[b]{\smash{{\SetFigFont{12}{14.4}{\rmdefault}{\mddefault}{\updefault}$j$}}}}
\put(3331,749){\makebox(0,0)[b]{\smash{{\SetFigFont{12}{14.4}{\rmdefault}{\mddefault}{\updefault}$j$}}}}
\put(2701,-1726){\makebox(0,0)[b]{\smash{{\SetFigFont{12}{14.4}{\rmdefault}{\mddefault}{\updefault}$j$}}}}
\put(3151,-1726){\makebox(0,0)[b]{\smash{{\SetFigFont{12}{14.4}{\rmdefault}{\mddefault}{\updefault}$j$}}}}
\put(3331,-1726){\makebox(0,0)[b]{\smash{{\SetFigFont{12}{14.4}{\rmdefault}{\mddefault}{\updefault}$k$}}}}
\put(3241,-286){\makebox(0,0)[b]{\smash{{\SetFigFont{12}{14.4}{\rmdefault}{\mddefault}{\updefault}$k$}}}}
\put(2611,-691){\makebox(0,0)[b]{\smash{{\SetFigFont{12}{14.4}{\rmdefault}{\mddefault}{\updefault}$i+j$}}}}
\put(1981,-286){\makebox(0,0)[b]{\smash{{\SetFigFont{12}{14.4}{\rmdefault}{\mddefault}{\updefault}$i+j$}}}}
\put(4231,749){\makebox(0,0)[b]{\smash{{\SetFigFont{12}{14.4}{\rmdefault}{\mddefault}{\updefault}$i$}}}}
\put(4861,-286){\makebox(0,0)[b]{\smash{{\SetFigFont{12}{14.4}{\rmdefault}{\mddefault}{\updefault}$i$}}}}
\put(4231,-691){\makebox(0,0)[b]{\smash{{\SetFigFont{12}{14.4}{\rmdefault}{\mddefault}{\updefault}$i$}}}}
\put(5491,-691){\makebox(0,0)[b]{\smash{{\SetFigFont{12}{14.4}{\rmdefault}{\mddefault}{\updefault}$i$}}}}
\put(4411,749){\makebox(0,0)[b]{\smash{{\SetFigFont{12}{14.4}{\rmdefault}{\mddefault}{\updefault}$j$}}}}
\put(5041,-286){\makebox(0,0)[b]{\smash{{\SetFigFont{12}{14.4}{\rmdefault}{\mddefault}{\updefault}$j$}}}}
\put(4411,-691){\makebox(0,0)[b]{\smash{{\SetFigFont{12}{14.4}{\rmdefault}{\mddefault}{\updefault}$j$}}}}
\put(5041,-1726){\makebox(0,0)[b]{\smash{{\SetFigFont{12}{14.4}{\rmdefault}{\mddefault}{\updefault}$j$}}}}
\put(5671,-691){\makebox(0,0)[b]{\smash{{\SetFigFont{12}{14.4}{\rmdefault}{\mddefault}{\updefault}$j$}}}}
\put(5491,-1726){\makebox(0,0)[b]{\smash{{\SetFigFont{12}{14.4}{\rmdefault}{\mddefault}{\updefault}$k$}}}}
\put(5671,-1726){\makebox(0,0)[b]{\smash{{\SetFigFont{12}{14.4}{\rmdefault}{\mddefault}{\updefault}$l$}}}}
\put(5491,749){\makebox(0,0)[b]{\smash{{\SetFigFont{12}{14.4}{\rmdefault}{\mddefault}{\updefault}$i$}}}}
\put(5671,749){\makebox(0,0)[b]{\smash{{\SetFigFont{12}{14.4}{\rmdefault}{\mddefault}{\updefault}$j$}}}}
\put(5491,-286){\makebox(0,0)[b]{\smash{{\SetFigFont{12}{14.4}{\rmdefault}{\mddefault}{\updefault}$k$}}}}
\put(5671,-286){\makebox(0,0)[b]{\smash{{\SetFigFont{12}{14.4}{\rmdefault}{\mddefault}{\updefault}$l$}}}}
\put(4861,-1726){\makebox(0,0)[b]{\smash{{\SetFigFont{12}{14.4}{\rmdefault}{\mddefault}{\updefault}$i$}}}}
\put(3511,-1186){\makebox(0,0)[lb]{\smash{{\SetFigFont{12}{14.4}{\rmdefault}{\mddefault}{\updefault}$\ul{i}+\ul{j}$}}}}
\put(3511,254){\makebox(0,0)[lb]{\smash{{\SetFigFont{12}{14.4}{\rmdefault}{\mddefault}{\updefault}$\ul{i}+\ul{k}$}}}}
\put(5851,254){\makebox(0,0)[lb]{\smash{{\SetFigFont{12}{14.4}{\rmdefault}{\mddefault}{\updefault}$\ul{i}+\ul{k}$}}}}
\put(5851,-1186){\makebox(0,0)[lb]{\smash{{\SetFigFont{12}{14.4}{\rmdefault}{\mddefault}{\updefault}$\ul{i}+\ul{k}$}}}}
\put(4141,-1726){\makebox(0,0)[b]{\smash{{\SetFigFont{12}{14.4}{\rmdefault}{\mddefault}{\updefault}$i+j$}}}}
\put(4771,-691){\makebox(0,0)[b]{\smash{{\SetFigFont{12}{14.4}{\rmdefault}{\mddefault}{\updefault}$i+j$}}}}
\put(4771,749){\makebox(0,0)[b]{\smash{{\SetFigFont{12}{14.4}{\rmdefault}{\mddefault}{\updefault}$i+j$}}}}
\put(4141,-286){\makebox(0,0)[b]{\smash{{\SetFigFont{12}{14.4}{\rmdefault}{\mddefault}{\updefault}$i+j$}}}}
\put(5086,749){\makebox(0,0)[b]{\smash{{\SetFigFont{12}{14.4}{\rmdefault}{\mddefault}{\updefault}$i$}}}}
\put(4456,-286){\makebox(0,0)[b]{\smash{{\SetFigFont{12}{14.4}{\rmdefault}{\mddefault}{\updefault}$i$}}}}
\put(4456,-1726){\makebox(0,0)[b]{\smash{{\SetFigFont{12}{14.4}{\rmdefault}{\mddefault}{\updefault}$i$}}}}
\end{picture}%
\end{center}

\noindent The chosen values simplify the computations greatly. Indeed, normally, there are three inequalities to check for each rule: hence, there should be 201 inequalities to check here. The first reduction comes from the fact that $F$ identifies $\tau$ and $\kappa$: there are 24 rules that can be dropped since, for each one, there is another rule that is sent to the same image. Thus there remains 43 rules and 129 inequalities to check.

Moreover, the rules of $\lz$ have some interesting symmetries that one can exploit: indeed, whenever $f\fl g$ is a rule of $\lz$, then $f^o\fl g^o$ is also a rule of $\lz$, where the duality $(\cdot)^o$ is the involution defined by:
$$
\mu^o=\delta, \quad \eta^o=\epsilon, \quad \tau^o=\tau, \quad \kappa^o=\kappa, \quad n^o=n, \quad (g\circ f)^o=f^o\circ g^o, \quad (f\tens g)^o=f^o\tens g^o.
$$

\noindent Another way to define this duality is by its action on diagrams: there, it is the top-down symmetry. Furthermore, the functor $F$ is compatible with this symmetry, in the sense that, for every arrow $f$, the functor $F$ sends $f^o$ onto $F(f)^o$, where the duality on $\Or$ is defined that way: $(f_*,f^*,[f])^o=(f^*,f_*,[f]^o)$, with $[f]^o(\vec{x},\vec{x}')=[f](\vec{x}',\vec{x})$. Note that this only have a meaning because the two sets $X$ and $Y$ are the same here (both equal to $\Nb^*$).

Thus, if some rule $f\fl g$ in $\lz$ satisfies $F(f)>F(g)$, then so does $f^o\fl g^o$. As a consequence, this reduces the number of rules to study: 18 of the remaining rules have a distinct dual, hence only 25 rules need to be studied (75 inequalities). Furthermore, when a rule $f\fl g$ is self-dual, the inequality $F(f)^*\geq F(g)^*$ holds if and only if $F(f)_*\geq F(g)_*$ holds: 8 of the remaining rules are in that case, which means there still are 67 inequalities from the former 201 to check. Computations do not rise any difficulty. For example, let us study the following (self-dual) rule:
\begin{center}
\begin{picture}(0,0)%
\includegraphics{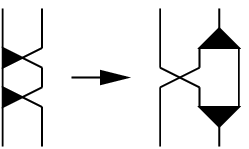}%
\end{picture}%
\setlength{\unitlength}{4144sp}%
\begingroup\makeatletter\ifx\SetFigFont\undefined%
\gdef\SetFigFont#1#2#3#4#5{%
  \reset@font\fontsize{#1}{#2pt}%
  \fontfamily{#3}\fontseries{#4}\fontshape{#5}%
  \selectfont}%
\fi\endgroup%
\begin{picture}(1104,654)(79,107)
\end{picture}%
\end{center}

\noindent One computes $\quad \begin{cases}
(\kappa\circ\kappa)_*(i,j)=(2i+j,i+j) \\
((1\tens\mu)\circ(\tau\tens 1)\circ(1\tens\delta))_*(i,j)=(i+j,i+j).
\end{cases}$ \\

\noindent Since $i$ and $j$ are non-zero natural numbers, the following inequality holds:
$$
(\kappa\circ\kappa)_*>((1\tens\mu)\circ(\tau\tens 1)\circ(1\tens\delta))_*.
$$

\noindent Then $\quad \begin{cases}
[\kappa\circ\kappa](i,j,k,l)=\ul{i}+\ul{i+j}+\ul{k}+\ul{k+l} \\
[(1\tens\mu)\circ(\tau\tens 1)\circ(1\tens\delta)](i,j,k,l)=2\ul{i}+\ul{j}+2\ul{k}+\ul{l}.
\end{cases}$ \\

\noindent Since $i$ and $j$ are non-zero natural numbers, the inequalities $i+j>i$ and $i+j>j$ always hold. Thus, by property of the multiset order on $[\Nb^*]$, the inequality $\ul{i+j}>\ul{i}+\ul{j}$ always holds. Similarly, so does $\ul{k+l}>\ul{k}+\ul{l}$. Finally, the multiset order on $[\Nb^*]$ is compatible with addition, yielding:
$$
[\kappa\circ\kappa]>[(1\tens\mu)\circ(\tau\tens 1)\circ(1\tens\delta)].
$$

\noindent The other rules are studied in a similar way [Guiraud 2004], which leads to the following result, proving that commutative equational theories can admit \emph{polygraphic} convergent presentations:

\begin{thm}\label{thm3}
The $3$-polygraph $\lz$ is a convergent presentation of the equational theory of $\zb$-vector spaces, with explicit resource management.
\end{thm}

\section*{Comments and future directions}
\emptysectionmark{Comments and future directions}

\noindent The study of ($3$-)polygraphs has been started by Albert Burroni and Yves Lafont, as an algebraic model for $3$-dimensional calculus on $2$-dimensional objects. Foundations were laid in [Lafont 1992], [Burroni 1993] and [Lafont 1995]. In [Lafont 2003], rewriting systems generated by $3$-polygraphs were considered and many known equational presentations are studied in order to be completed into convergent rewriting systems (or, at least, rewriting systems with the unique normal form property). Discussions with Albert Burroni, Yves Lafont and Philippe Malbos have been essential in order to achieve the results presented here. Comments from the referee were of great help to make this paper clearer.

There exist many research paths concerning polygraph. The first one is about confluence: as mentionned earlier, there exist theoretical issues with critical pairs of $3$-polygraphs; exploration and classification are mandatory in order to achieve some automated completion procedure for these objects. Such a tool (which implementation in Caml has already started) would be very useful since, starting from an equational theory, one could use the constructions described in section \ref{explicit} in order to obtain a $3$-polygraph; then a completion procedure could be applied to correct termination and confluence issues. Suggested by Pierre Lescanne, other usual techniques for building reduction orders in term rewriting could also be examined, in order to see if they could also be adaptated to polygraphs. Among the most useful results to be studied are the ones concerning path orders, see [Baader Nipkow 1998], and dependency pairs, see [Arts Giesl 2000].

A second theme to be explored is the study of higher-dimensional polygraphs. For an example of application, $4$-polygraphs provide a categorical framework for proof transformations in the \emph{calculus of structures} [Guglielmi Stra\ss{}burger 2001]. Such an approach could yield results such as proof decompositions or normal forms, given by a convergent $4$-polygraph. At least, it suggests that formulas are $2$-dimensional objects, proofs are $3$-dimensional and computation on them (such as cut elimination) lives in dimension $4$. This point of view is conjectured to yield a new class of objects describing formal proofs, giving a different, categorical and geometrical way to approach proof theory.

Theoretical studies can also be directed at pursuing the synthesis started in [Guiraud 2004] on rewriting systems: one of the main goals is to have a framework where one can compare two rewriting systems, regardless of the algebraic structure of their terms. The reduction space associated to each rewriting system is an algebraico-geometric object (a cubical object in some category of algebras) and one could use the underlying cubical sets of these objects to compare rewriting systems, geometrically. Notions of (co)fibrations from Quillen model categories - see [Hovey 1999] - theory could be useful for a better understanding of results such as the ones of section \ref{explicit}; since many rewriting systems are special cases of polygraphs, this study will start with the construction of homotopical tools for these objects.

Still another question is the following: is there some $n$ for which there exists a finite $n$-polygraph yielding a calculus with both explicit substitutions and explicit resource management for the $\lambda$-calculus. When $n=3$, the answer seems to be negative, since theoretical results deny the existence of any non-trivial product category that is both cartesian (for resource management) and sovereign (for substitutions). An equational description of the structure of closed category (such as the one Albert Burroni has given for cartesian categories) should be the first step of this work. Another possibility is to use a $3$-dimensional interpretation of proofs, together with the links between $\lambda$-terms and proofs.

Finally, $3$-polygraphs have the interesting property to modelize computational circuits. Indeed, both classical and quantum algorithms accept representations as circuits which are, albeit not in their usual presentation, genuine operators of a $3$-polygraph. Furthermore, equational presentations are known for both kinds of circuits. Questions that can be studied with this point of view concern the existence of convergent $3$-polygraphs for classical or quantum circuits, thus leading to canonical representations of programs. One can take a look at [Kitaev Shen Vyalyi 2002] for more information on circuits and [Lafont 2003] for their links with polygraphs.

\section*{References}
\emptysectionmark{References}

\noindent \textsc{T. Arts} and \textsc{J. Giesl}, \emph{Termination of term rewriting using dependency pairs}. Theoretical Computer Science 236, 133-178, 2000.

\smallskip
\noindent \textsc{F. Baader} and \textsc{T. Nipkow}, \emph{Term rewriting and all that}. Cambridge University Press, 1998.

\smallskip
\noindent \textsc{A. Burroni}, \emph{Higher-dimensional word problems with applications to equational logic}. Theoretical Computer Science 115(1), 46-62, 1993.

\smallskip
\noindent \textsc{A. Guglielmi and} \textsc{L. Stra\ss{}burger}, \emph{Non-commutativity and MELL in the calculus of structures}. Lecture Notes in Computer Science 2142, 54-68, 2001.

\smallskip
\noindent \textsc{Y. Guiraud}, \emph{Présentations d'opérades et systèmes de réécriture}. Thèse de doctorat, Montpellier, 2004.

\smallskip
\noindent \textsc{M. Hovey}, \emph{Model categories}. Mathematical Surveys and Monographs 63, 1999.

\smallskip
\noindent \textsc{A. Kitaev}, \textsc{A. Shen} and \textsc{M. Vyalyi}, \emph{Classical and quantum computation}. Graduate Studies in Mathematics 47, 2002.

\smallskip
\noindent \textsc{Y. Lafont}, \emph{Penrose diagrams and 2-dimensional rewriting}. London Mathematical Society Lecture Notes Series 177, 191-201, 1992.

\smallskip
\noindent ---, \emph{Equational reasoning with 2-dimensional diagrams}. Lecture Notes in Computer Science 909, 170-195, 1995.

\smallskip
\noindent ---, \emph{Towards an algebraic theory of boolean circuits}. Journal of Pure and Applied Algebra 184, 257-310, 2003.

\smallskip
\noindent \textsc{S. MacLane}, \emph{Categorical algebra}. Bulletin of the American Mathematical Society 71, 40-106, 1965.

\smallskip
\noindent ---, \emph{Categories for the working mathematician}. Springer, 1998.

\end{document}